\documentclass{article}%
\usepackage{amsmath}
\usepackage{amsfonts}
\usepackage{amssymb}
\usepackage{amsthm}
\usepackage{graphicx}
\usepackage{subcaption}
\usepackage{color}
\usepackage{hyperref}
\usepackage{xcolor}
\usepackage{float}
\definecolor{mygreen}{rgb}{0,0.5,0}

\newtheorem{theorem}{Theorem}

\newtheorem{proposition}[theorem]{Proposition}
\newtheorem{remark}[theorem]{Remark}

\title{
From Low-Dimensional Chaos to High-Dimensional Dynamics \\
in Self-Attention Systems
}

\author{
J. D. Mireles James 
 \thanks{J.D.M.J. partially supported by NSF grant DMS - 2307987.}
\thanks{Florida Atlantic University, Department of Mathematical Sciences,
Science Building, Room 234, 777 Glades Road,
Boca Raton, Florida, 33431, USA, {\tt jmirelesjames@fau.edu}.}
}

\date{November 22, 2026}

\begin{document}

\maketitle

\begin{abstract}
Continuous-time particle models of self-attention are often studied under symmetry assumptions that 
yield a gradient structure and constrain the dynamics to convergence toward equilibria. We investigate 
the behavior of these systems beyond the gradient setting, considering autonomous projected self-attention 
models with constant query, key, and value matrices and particles constrained to an ellipse. For three particles, 
the equations define a vector field on the three-dimensional torus. We derive explicit synchronized periodic 
solutions and compute their Floquet multipliers, establishing the existence of non-gradient dynamics. 
Numerical explorations using bifurcation diagrams, Poincaré sections, periodic-orbit continuation, and Lyapunov exponents reveal attracting periodic orbits, invariant-torus-like behavior, and strong evidence of chaos. We then 
construct higher-dimensional systems by coupling these low-dimensional building blocks through 
masked and multi-head attention mechanisms. We discuss the persistence of normally hyperbolic invariant 
manifolds under weak coupling and investigate how the coupling architecture influences the resulting dynamics. 
Numerical experiments reveal the emergence of additional expanding directions, including a 60-particle system 
in which 26 of the 40 computed Lyapunov exponents exceed 0.01. These results demonstrate the richness 
of the dynamical behavior possible in simplified self-attention systems and motivate further investigation 
of the invariant structures organizing their high-dimensional dynamics.
\end{abstract}



\section{Introduction}
\label{sec:introduction}

The 2017 paper \emph{Attention Is All You Need} by Vaswani et al. \cite{VaswaniEtAl2017}
introduced a novel network architecture 
called a \textit{transformer} and spurred rapid advances in
text-based generative artificial intelligence (AI). In the years since, 
transformers have become ubiquitous in Large Language Models (LLMs) 
including OpenAI's ChatGPT,  Anthropic's Claude, 
and Meta's Llama \cite{OpenAI2023GPT4,Lindsey2025,TouvronEtAl2023} to 
name only a few.  Indeed, \emph{GPT} stands for Generative Pre-trained Transformer
\cite{BrownEtAl2020}.

Very roughly speaking, the attention mechanism (formalized in Equation 
\eqref{eq:intro-attention-weights} below) treats a language token as a particle 
on a compact manifold and a text prompt as a cloud of such particles.  
The transformer is a layered network which reorganizes the point cloud, so that 
an appropriate average of the output state provides the most likely next token. 
As the cloud of particles passes through the network, each particle nudges --- and is 
nudged --- by the others.  The strength of the adjustment is determined by the 
current state of the whole cloud.  Projection onto the constraint manifold prevents 
unbounded growth of particle trajectories.

A concise mathematical description of this operation
replaces network depth, in an appropriate limit, by continuous time.
This results in a system of differential equations on
the constraint manifold (see Equation \eqref{eq:intro-particle-equations}
below), where field interactions between particles encode 
prompt context.  That is, the same language token
follows a different path, and shapes the 
 trajectories of its neighbors differently, depending on the relative positions 
 of the particles.  

The particle-system viewpoint of attention has roots in the work of 
Lu et al.~\cite{LuEtAl2019} and in the measure-theoretic
treatments of Vuckovic et al.~\cite{VuckovicEtAl2020} and also 
Sander et al.~\cite{SanderEtAl2022}. 
For more discussion from a decidedly mathematical physics perspective, we
refer to the work of Geshkovski, Letrouit, Polyanskiy, and Rigollet
\cite{GeshkovskiEtAl2023,GeshkovskiLetrouitPolyanskiyRigollet2025}, whose
approach we follow below, and to Polyanskiy, Rigollet, and Yao
\cite{PolyanskiyRigolletYao2025} for the case of the circle.

To formalize the discussion, let $d\geq 2$ be the dimension of the token embedding 
space, and $W\in\mathbb R^{d\times d}$ be a symmetric, positive-definite 
matrix.  Define the ellipsoid
\[
\mathcal E_W=\{x\in\mathbb R^d:x^\top W x=1\}.
\]
and the projection operator $\Pi^w_x \colon \mathbb{R}^d \to T_x \mathcal{E}_W$
by 
\[
\Pi^W_xv=v-x(x^\top Wv).
\]

Let $N$ denote the width of the context window, that is, 
the number of language tokens in the prompt (i.e. the number of
particles in the system). We write 
$x_1,\ldots,x_N \in \mathcal E_W$ to denote the individual
tokens, and $X=(x_1,\ldots,x_N)\in\mathcal E_W^N$ to represent
the entire prompt.  That is, $X$ represents a state of the 
self-attention transformer.  

Let $Q,K:[0,T]\to\mathbb R^{m\times d}$ and
$V:[0,T]\to\mathbb R^{d\times d}$ be smooth, matrix-valued functions.
These are referred to as the query, key, and value matrices, respectively.
(See Remark \ref{rem:Physics} below for the meaning of these terms
in machine learning applications).
For $\beta\geq 0$, define the attention coefficients, 
or attention weights, by 
\begin{equation}\label{eq:intro-attention-weights}
\alpha_{ij}(t,X)=
\frac{\exp\!\bigl(\beta\langle Q(t)x_i,K(t)x_j\rangle\bigr)}
{\displaystyle\sum_{k=1}^N
\exp\!\bigl(\beta\langle Q(t)x_i,K(t)x_k\rangle\bigr)},
\qquad 1\leq i,j\leq N.
\end{equation}
The denominator is referred to as the normalization of the attention weight, and 
the parameter $1/\beta$ is referred to as the temperature of the system. 
The equations of motion for the  
continuous-time, single-head, self-attention particle
system are given by the following projected vector field 
\begin{equation}\label{eq:intro-particle-equations}
\dot x_i(t)=
\Pi^W_{x_i(t)}
\left(
\sum_{j=1}^N
\alpha_{ij}(t,X(t))V(t)x_j(t)
\right),
\qquad 1\leq i\leq N.
\end{equation}
These provide the equations of motion for the 
particle system on the ellipsoid.  The output of the transformer is
an appropriate average of the state of the system after some 
specified time $T > 0$.

\begin{remark}[Interpretation of the attention mechanism] \label{rem:Physics}
These choices, and the accompanying terminology, are more thoroughly
justified and contextualized in
\cite{GeshkovskiEtAl2023,GeshkovskiLetrouitPolyanskiyRigollet2025}.
The interested reader is referred also to
\cite{LuEtAl2019,RodriguezAbellaSilvestreTabuadaAsymptotic,Burger2025,CastinEtAl2025}
and the references therein.
Nevertheless, we formulate a few brief clarifying remarks.

The query, key, and value matrices have distinct roles in the attention
mechanism of Vaswani et al.~\cite{VaswaniEtAl2017}.  The query $Q(t)x_i$
represents the information sought by particle $i$, while the key $K(t)x_j$
represents the information made available by particle $j$.  Their inner
product measures the compatibility of these quantities, determining the
relative strength of the particle interaction.  The value $V(t)x_j$
represents the information transmitted by particle $j$ to particle $i$.
Thus, the attention mechanism is a weighted average of the token values,
with weights determined by query-key compatibility.

The softmax normalization is common in statistical mechanics and information
theory.  For fixed $i$, writing $s_{ij}=\langle Q(t)x_i,K(t)x_j\rangle$, the
attention weights form a Gibbs probability distribution,
\[
\alpha_{ij}=\frac{e^{\beta s_{ij}}}{\sum_k e^{\beta s_{ik}}}.
\]
For $\beta>0$ (that is, temperature less than infinity), the Gibbs
distribution uniquely maximizes
\[
\sum_j p_j s_{ij}+\beta^{-1}H(p),
\]
over the probability simplex, where $H(p)=-\sum_j p_j\log p_j$ is Shannon
entropy.  (This functional is the negative free energy, with energies
$-s_{ij}$ and temperature $1/\beta$.)  Thus, $\beta$ controls competition
between selecting strongly compatible particles and distributing attention
more uniformly.  As $\beta\to0$, the weights approach $1/N$, while as
$\beta\to\infty$, attention concentrates on the particles with maximal
compatibility.  The weights are positive and sum to one, and they depend
smoothly on the particle configuration, since the exponential is smooth and
the denominator is positive.
\end{remark}

\begin{remark}[Model fidelity] \label{rem:realModels}
The system given by Equation \eqref{eq:intro-particle-equations} 
is an idealization of the architecture introduced in \cite{VaswaniEtAl2017}, and 
of those used in practical applications.  The constraint to the ellipsoid
plays the role of layer normalization, but practical transformers differ
from the model in at least three ways.  First, they use multiple attention
heads, each with its own query, key, and value matrices; in continuous time
this corresponds to a sum of single-head terms.  Second, attention layers
are often combined with feedforward layers, adding a velocity field which acts
on each particle separately.  In an appropriate continuous-depth limit, 
such layers lead to neural-ODE networks \cite{ChenEtAl2018}.
Third, attention is often masked, so that each particle interacts with only
a subset of the others.  Causal masking in autoregressive language models,
where each token attends only to the tokens preceding it, is one example
\cite{Karagodin2024,CastinEtAl2025}; masks encoding graph or type structure,
as in graph and heterogeneous attention architectures
\cite{VelickovicEtAl2018,HuEtAl2020}, are another.  

This said, the system of Equation \eqref{eq:intro-particle-equations} is 
nevertheless a \textit{building block} for more realistic systems like the
multi-head attention actually employed in \cite{VaswaniEtAl2017}.
In the present work (see  Section~\ref{sec:bigTransformers})
we use masked multi-head attention to assemble 
more realistic systems out of the building blocks 
considered in Section \ref{sec:torus}.
\end{remark}

\bigskip

An important qualitative result informing our 
understanding of the emergence of consensus, 
clustering, and metastability in self-attention systems
is that, under appropriate symmetry 
restrictions on the query, key, and value matrices, these systems 
have gradient structure. 
Such structure severely constrains the dynamics: in an analytic gradient
system on a compact manifold, every orbit converges, in both forward and
backward time, to an equilibrium (not necessarily the same one).
 The theorem below follows, by a
linear change of variables, from the spherical gradient-flow construction
of Geshkovski, Letrouit, Polyanskiy, and Rigollet
\cite[Section~1.2.4]{GeshkovskiLetrouitPolyanskiyRigollet2025}.

\begin{theorem}[Gradient structure of symmetric self-attention]
\label{thm:intro-gradient}
Suppose that $\beta>0$ and that $Q$, $K$, and $V$ are constant, with
\[
 P:=Q^\top K=WV,\qquad P=P^\top.
\]
Define
\[
 Z_i(X)=\sum_{j=1}^N\exp(\beta x_i^\top P x_j),\qquad
 \mathcal H_\beta(X)=\frac{1}{2\beta}
 \sum_{i,j=1}^N\exp(\beta x_i^\top P x_j),
\]
and equip $\mathcal E_W^N$ with the Riemannian metric
\[
 g_X(\xi,\eta)=\sum_{i=1}^N Z_i(X)\,\xi_i^\top W\eta_i.
\]
Then \eqref{eq:intro-particle-equations} is the gradient-ascent flow of
$\mathcal H_\beta$ with respect to $g$. Along every solution,
\[
 \frac{d}{dt}\mathcal H_\beta(X(t))=
 \sum_{i=1}^N Z_i(X(t))\,\dot x_i(t)^\top W\dot x_i(t)\geq0,
\]
with equality if and only if $X(t)$ is an equilibrium. In particular,
there are no nonconstant periodic solutions, and every trajectory
converges to an equilibrium in forward and backward time.
\end{theorem}

The last assertion follows from analyticity, compactness, and the
\L{}ojasiewicz gradient inequality. The theorem and its generalizations
illuminate important aspects of clustering and the long-time organization
of token configurations. For related results and complementary viewpoints,
we refer to \cite{GeshkovskiEtAl2023,GeshkovskiEtAl2024,
GeshkovskiLetrouitPolyanskiyRigollet2025,
RodriguezAbellaSilvestreTabuadaAsymptotic,CastinEtAl2025,KuehnYoon2026}.

\bigskip

While Theorem \ref{thm:intro-gradient} provides a powerful 
description of the qualitative dynamics, the algebraic hypotheses are in fact
quite restrictive. Typical applications to machine learning
impose no such constraints on 
 the query, key, and value matrices during training~\cite{VaswaniEtAl2017}, 
 and since $Q^\top K=WV$ and the symmetry of $Q^\top K$ 
are nongeneric in the set of real matrices, these conditions 
should be expected to fail in general.  
We refer to Tomihari and
Karakida~\cite{TomihariKarakida2025}
for another perspective on self-attention
in the absence of an energy function.

The goal of the present work is to explore how far transformer
dynamics can depart from this gradient picture.  To this end,
we consider the elementary case of three particles on an ellipse
in $\mathbb{R}^2$, with constant but otherwise unrestricted attention
matrices.  The right-hand side of
Equation~\eqref{eq:intro-particle-equations} defines a vector field on
$\mathbb R^6$ that is tangent to the embedded torus $\mathcal E_W^3$.
In angular coordinates, that is, on the universal cover $\mathbb R^3$,
it becomes a system of three unprojected ordinary differential equations,
$2\pi$-periodic in each angle variable, and
the dynamics of three-dimensional vector fields
are much easier to visualize. 

We begin our study by giving conditions under which this
three-dimensional example system possesses nonconstant periodic orbits.
For such parameters the system admits no strict Lyapunov function on the
torus and, in particular, is not the gradient flow of any single-valued
energy with respect to any Riemannian metric. 
Consequently, such parameters necessarily violate the hypotheses of
Theorem~\ref{thm:intro-gradient}.

Next, we present numerical explorations of the three-particle system using
classical tools from computational dynamics such as
Poincar\'e sections, Lyapunov exponents, periodic-orbit continuation, and
bifurcation diagrams. The computations reveal periodic windows, changes
of stability of periodic orbits, and provide compelling numerical evidence 
for the existence of chaotic attractors. 

Finally, we show how elementary low-dimensional attention 
systems can be assembled into higher-dimensional models using masked
multi-head attention.  That is, we construct examples where,
at zero coupling, we have a product of independently evolving
systems whose dynamics we choose.  Increasing 
the coupling parameter introduces interactions between the constituent subsystems, 
allowing us to build large networks with complex, non-gradient dynamics.
We study these networks using standard ideas from perturbation theory, as well as numerical 
experiments.  These constructions are motivated by established
attention architectures involving restricted and type-dependent
interactions \cite{VelickovicEtAl2018,HuEtAl2020,CastinEtAl2025}.
See Section \ref{sec:structureApplications} for more detailed discussion.

A model related to the one derived in Section~\ref{sec:torus} was recently
studied by Li et al.~\cite{LiEtAl2026LinearTransformers}, who consider
particles constrained to the circle with a linear attention kernel in place
of the exponential softmax.  They show that the resulting dynamics admit a
generalized Kuramoto formulation with second-harmonic coupling, and use
classical methods for coupled oscillators to investigate clustering,
periodic oscillations, and bifurcations, primarily in the mean-field limit.
Here we retain the softmax attention mechanism, work with finitely many
particles, and construct specific examples with complicated, non-gradient 
dynamics.

\medskip

\noindent \textbf{Code and data availability:} 
The MATLAB implementations used in this paper, together with archived numerical data and scripts for reproducing the figures and tables, are publicly available through the accompanying project webpage
\cite{MirelesJames2026AttentionDynamicsWebsite}.

\section{A three particle model on the torus}\label{sec:torus}

Throughout this section, we consider three particles $x_1, x_2, x_3 \in \mathbb{R}^2$ and
constant query, key, and value matrices. The particles are constrained to 
the surface of an ellipse defined as follows.  Let $W\in\mathbb R^{2\times2}$ 
be symmetric and positive definite, and write
\[
\mathcal E_W=\{x\in\mathbb R^2:x^\top Wx=1\},\qquad S=W^{-1/2}.
\]
The map $u\mapsto Su$ takes the unit circle diffeomorphically onto
$\mathcal E_W$. We introduce
\[
u(\theta)=\begin{pmatrix}\cos\theta\\\sin\theta\end{pmatrix},
\qquad
v(\theta)=u'(\theta)=\begin{pmatrix}-\sin\theta\\\cos\theta\end{pmatrix},
\]
and parametrize the particle positions by
\[
x_i=Su(\theta_i),\qquad \theta_i\in\mathbb R/2\pi\mathbb Z,
\qquad i=1,2,3.
\]
In particular, $\dot x_i=Sv(\theta_i)\dot\theta_i$, where the dot denotes
differentiation with respect to time. We identify the configuration space
$\mathcal E_W^3$ with the three-dimensional torus
$\mathbb T^3=(\mathbb R/2\pi\mathbb Z)^3$.

The matrices $Q$, $K$, and $V$ enter the angular equations through the
effective $2\times2$ matrices
\begin{equation}\label{eq:torus-effective-matrices}
A=W^{1/2}VW^{-1/2},\qquad
B=W^{-1/2}Q^\top KW^{-1/2}.
\end{equation}
For convenience, put $u_i=u(\theta_i)$ and $v_i=v(\theta_i)$. The
attention coefficients become
\begin{equation}\label{eq:torus-weights}
\alpha_{ij}(\theta)
=\frac{\exp\!\bigl(\beta u_i^\top Bu_j\bigr)}
{\displaystyle\sum_{k=1}^3\exp\!\bigl(\beta u_i^\top Bu_k\bigr)},
\qquad 1\leq i,j\leq3.
\end{equation}
To obtain the angular equations, observe that $Sv_i$ has unit length
in the $W$-metric, since $(Sv_i)^\top W(Sv_i)=1$. Moreover, for
$w\in\mathbb R^2$, the tangent projection
$\Pi^W_{x_i}w=w-x_i(x_i^\top Ww)$ satisfies
\[
(Sv_i)^\top W\Pi^W_{x_i}w
=(Sv_i)^\top Ww=v_i^\top W^{1/2}w.
\]
Taking the $W$-inner product of \eqref{eq:intro-particle-equations}
with $Sv_i$, and using \eqref{eq:torus-effective-matrices}, gives the
autonomous system
\begin{equation}\label{eq:torus-angular-equations}
\dot\theta_i=F_i(\theta)
:=v_i^\top A\sum_{j=1}^3\alpha_{ij}(\theta)u_j
=\frac{\displaystyle\sum_{j=1}^3
 e^{\beta u_i^\top Bu_j}\,v_i^\top Au_j}
 {\displaystyle\sum_{k=1}^3e^{\beta u_i^\top Bu_k}},
\qquad i=1,2,3.
\end{equation}

For later numerical work, it is useful to write the right-hand side
completely in terms of trigonometric functions and the entries of $A$
and $B$. Let
\[
A=\begin{pmatrix}a_{11}&a_{12}\\a_{21}&a_{22}\end{pmatrix},
\qquad
B=\begin{pmatrix}b_{11}&b_{12}\\b_{21}&b_{22}\end{pmatrix},
\qquad c_i=\cos\theta_i,\quad s_i=\sin\theta_i.
\]
Then
\[
u_i^\top Bu_j
=b_{11}c_ic_j+b_{12}c_is_j+b_{21}s_ic_j+b_{22}s_is_j,
\]
whereas
\[
v_i^\top Au_j
=(a_{21}c_i-a_{11}s_i)c_j
 +(a_{22}c_i-a_{12}s_i)s_j.
\]
Substitution in \eqref{eq:torus-angular-equations} yields the explicit
three-dimensional vector field
\begin{equation}\label{eq:torus-explicit-trigonometric}
\dot\theta_i=
\frac{\displaystyle\sum_{j=1}^3
\exp\!\left(\beta\left[b_{11}c_ic_j+b_{12}c_is_j
+b_{21}s_ic_j+b_{22}s_is_j\right]\right)
\left[(a_{21}c_i-a_{11}s_i)c_j
+(a_{22}c_i-a_{12}s_i)s_j\right]}
{\displaystyle\sum_{j=1}^3
\exp\!\left(\beta\left[b_{11}c_ic_j+b_{12}c_is_j
+b_{21}s_ic_j+b_{22}s_is_j\right]\right)},
\qquad i=1,2,3.
\end{equation}
Thus $\dot\theta=F(\theta)$ defines a real-analytic vector field on
$\mathbb T^3$. Its lift to $\mathbb R^3$ is $2\pi$-periodic in each
coordinate, and the denominators in
\eqref{eq:torus-explicit-trigonometric} are everywhere strictly
positive. In particular, the flow exists for all time on the compact
phase space $\mathbb T^3$.

The effective matrices $A$ and $B$ can be prescribed independently.
Indeed, for arbitrary real $2\times2$ matrices $A$ and $B$, one can
choose the query/key dimension $m=2$ and take
\[
V=W^{-1/2}AW^{1/2},\qquad Q=I_2,\qquad
K=W^{1/2}BW^{1/2}.
\]
Consequently, when the attention matrices are unrestricted, the ellipse
can be removed by a linear change of coordinates without restricting
the class of angular vector fields
\eqref{eq:torus-angular-equations}. Note also that the gradient
hypotheses $Q^\top K=WV$ and $(Q^\top K)^\top=Q^\top K$ from
Theorem~\ref{thm:intro-gradient} reduce in these coordinates to
$A=B=B^\top$.
 
Because all three particles use the same attention matrices, the vector
field is equivariant under permutations of the particles. In particular,
the pair-coincidence subtori
$\{\theta_i=\theta_j\pmod{2\pi}\}$ and the fully synchronized circle
\[
\Delta=\{(\varphi,\varphi,\varphi):\varphi\in\mathbb T\}
\]
are invariant. On $\Delta$, all attention weights are equal to $1/3$,
and the equations reduce to
\begin{equation}\label{eq:torus-synchrony}
\dot\varphi=v(\varphi)^\top Au(\varphi)
=a_{21}\cos^2\varphi-a_{12}\sin^2\varphi
+(a_{22}-a_{11})\sin\varphi\cos\varphi,
\end{equation}
which is independent of $B$ and $\beta$. The synchronized dynamics
provides an explicit family of nongradient examples, as the following 
proposition shows.

\begin{proposition}[Existence and stability of a synchronized periodic orbit]
\label{prop:torus-rotating-wave}
Consider the three-particle system on the unit circle,
written in angular coordinates as
\[
\dot\theta_i
=
\sum_{j=1}^3
\alpha_{ij}(\theta)
v(\theta_i)^T A u(\theta_j),
\qquad 1\leq i\leq 3,
\]
where
\[
u(\theta)=
\begin{pmatrix}\cos\theta\\\sin\theta\end{pmatrix},
\qquad
v(\theta)=u'(\theta),
\]
and
\[
\alpha_{ij}(\theta)
=
\frac{\exp\bigl(\beta u(\theta_i)^T B u(\theta_j)\bigr)}
{\sum_{k=1}^3
\exp\bigl(\beta u(\theta_i)^T B u(\theta_k)\bigr)}.
\]
Suppose that $A,B\in\mathbb R^{2\times 2}$ are constant,
$\beta\geq 0$, and that $A$ has complex conjugate eigenvalues
$a\pm i\kappa$, where
\[
a=\frac12\operatorname{tr}(A),
\qquad
\kappa=
\sqrt{\det(A)-\frac14\operatorname{tr}(A)^2}>0.
\]
Then the synchronized diagonal
\[
\mathcal D
=
\{(\varphi,\varphi,\varphi):
\varphi\in\mathbb T\}
\]
is a periodic orbit with period
\[
T=\frac{2\pi}{\kappa}.
\]
Its monodromy matrix is
\[
M=
\frac13\mathbf 1\mathbf 1^T
+
\sigma\left(I_3-\frac13\mathbf 1\mathbf 1^T\right),
\qquad
\sigma=e^{-aT},
\]
where $\mathbf 1=(1,1,1)^T$.
Equivalently,
\[
M=\frac13
\begin{pmatrix}
1+2\sigma&1-\sigma&1-\sigma\\
1-\sigma&1+2\sigma&1-\sigma\\
1-\sigma&1-\sigma&1+2\sigma
\end{pmatrix}.
\]
In particular, the Floquet multipliers are
\[
1,\qquad \sigma,\qquad \sigma.
\]

Moreover, for any $\varphi_0\in\mathbb T$, the Poincar\'e
return map on the section
$\Sigma=\{\theta_1=\varphi_0\}$ has a fixed point
corresponding to the synchronized orbit.
In local coordinates centered at this fixed point,
its derivative is
\[
DP(0,0)=\sigma I_2.
\]
Consequently, the synchronized periodic orbit is
hyperbolic and orbitally attracting if
$\operatorname{tr}(A)>0$, and hyperbolic and
orbitally repelling if $\operatorname{tr}(A)<0$.
If $\operatorname{tr}(A)=0$, the orbit is nonhyperbolic.
These conclusions are independent of $B$ and $\beta$.
\end{proposition}

\begin{proof}
At synchrony, $\theta_1=\theta_2=\theta_3=\varphi$,
all attention weights are equal to $1/3$. Hence
the synchronized diagonal is invariant, and the
dynamics restricted to it are governed by
\[
\dot\varphi=g(\varphi),
\qquad
g(\varphi)=v(\varphi)^TAu(\varphi).
\]
Since $A$ has no real eigenvalues, $g(\varphi)$
cannot vanish. Indeed, $g(\varphi)=0$ would imply
that $Au(\varphi)$ is parallel to $u(\varphi)$,
giving a real eigenvector of $A$. Thus the
synchronous motion is periodic.

To compute its period, observe that the equation
for $u(\varphi(t))$ can be written as
\[
\frac{d}{dt}u
=
(I_2-uu^T)Au.
\]
Consequently,
\[
u(\varphi(t))
=
\frac{e^{At}u(\varphi_0)}
{\|e^{At}u(\varphi_0)\|}.
\]
Writing $A=aI_2+C$, where $C^2=-\kappa^2I_2$,
we obtain
\[
e^{At}
=
e^{at}
\left(
\cos(\kappa t)I_2+
\frac{\sin(\kappa t)}{\kappa}C
\right).
\]
Since $A$ has no real eigenvectors, the first
positive time at which $e^{At}u(\varphi_0)$ is
a positive multiple of $u(\varphi_0)$ is
$T=2\pi/\kappa$. This is the period of the
synchronized orbit.

We next compute the variational equation along
the synchronized solution. Define
\[
d(\varphi)=u(\varphi)^TAu(\varphi),
\qquad
c(\varphi)=v(\varphi)^TAv(\varphi).
\]
At synchrony, the derivative of the attention
weights does not contribute to the linearization.
Indeed, the attention weights satisfy
$\sum_j\alpha_{ij}=1$, and all the vectors
being averaged coincide. Therefore the
corresponding first-order terms cancel.

Using $u'=v$ and $v'=-u$, direct differentiation
gives
\[
\frac{\partial F_i}{\partial\theta_k}
=
-d(\varphi)\delta_{ik}
+\frac13c(\varphi),
\]
where $F_i$ denotes the right-hand side of
the $i$th angular equation and $\delta_{ik}$
is the Kronecker delta. Thus the variational
equation is
\[
\dot\xi=L(\varphi(t))\xi,
\qquad
L(\varphi)
=
-d(\varphi)I_3
+\frac13c(\varphi)\mathbf 1\mathbf 1^T.
\]

Introduce the complementary projections
\[
P_{\parallel}
=
\frac13\mathbf 1\mathbf 1^T,
\qquad
P_{\perp}=I_3-P_{\parallel}.
\]
Since
\[
g'(\varphi)
=
-d(\varphi)+c(\varphi),
\]
the Jacobian is
\[
L(\varphi)
=
g'(\varphi)P_{\parallel}
-d(\varphi)P_{\perp}.
\]
The projections are constant and commute.
Consequently, the fundamental matrix of the
variational equation, normalized by
$\Phi(0)=I_3$, is
\[
\Phi(t)
=
\exp\left(\int_0^t g'(\varphi(s))\,ds\right)
P_{\parallel}
+
\exp\left(-\int_0^t d(\varphi(s))\,ds\right)
P_{\perp}.
\]

Consider these integrals over one period.
Since $\dot\varphi=g(\varphi)$ and $g$ never
vanishes,
\[
\frac{d}{dt}\log|g(\varphi(t))|
=
g'(\varphi(t)).
\]
The periodicity of $g$ therefore gives
\[
\int_0^T g'(\varphi(t))\,dt=0.
\]
On the other hand, because $u(\varphi)$ and
$v(\varphi)$ form an orthonormal basis of
$\mathbb R^2$,
\[
d(\varphi)+c(\varphi)=\operatorname{tr}(A)=2a.
\]
Combining this identity with
$g'=c-d$ gives
\[
d(\varphi)=a-\frac12g'(\varphi).
\]
It follows that
\[
\int_0^T d(\varphi(t))\,dt=aT.
\]
Hence the monodromy matrix is
\[
M=\Phi(T)
=
P_{\parallel}+e^{-aT}P_{\perp},
\]
which is precisely the stated formula.
The vector $\mathbf 1$ is an eigenvector
with eigenvalue $1$, while every vector
in the two-dimensional subspace
\[
\xi_1+\xi_2+\xi_3=0
\]
is an eigenvector with eigenvalue
$\sigma=e^{-aT}$. This proves the formula
for the Floquet multipliers.

It remains to compute the derivative of the
Poincar\'e return map. Fix $\varphi_0$ and
consider the transverse section
\[
\Sigma=\{\theta_1=\varphi_0\pmod{2\pi}\},
\]
with local coordinates
\[
(x,y)\longmapsto
(\varphi_0,\varphi_0+x,\varphi_0+y).
\]
Let $\Psi^t$ denote the flow of the full
three-particle system. Since
$g(\varphi_0)\ne0$, the first return time
$\tau(x,y)$ near the synchronized orbit is
well-defined and smooth, with
$\tau(0,0)=T$. Define
\[
P(x,y)
=
\pi_{2,3}
\Psi^{\tau(x,y)}
(\varphi_0,\varphi_0+x,\varphi_0+y)
-(\varphi_0,\varphi_0)
\pmod{2\pi},
\]
where $\pi_{2,3}$ denotes projection onto
the second and third angular coordinates.
Clearly $P(0,0)=(0,0)$.

At the synchronized crossing, the vector
field is $g(\varphi_0)\mathbf 1$.
Differentiating the return condition
$\theta_1(\tau(x,y))=\varphi_0\pmod{2\pi}$
gives the standard correction for variations
of the return time. Therefore
\[
DP(0,0)=
\pi_{2,3}
(I_3-\mathbf 1e_1^T)M
\begin{pmatrix}
0&0\\
1&0\\
0&1
\end{pmatrix},
\]
where $e_1=(1,0,0)^T$.
Substituting
$M=P_{\parallel}+\sigma P_{\perp}$,
and observing that
$(I_3-\mathbf 1e_1^T)P_{\parallel}=0$,
we obtain
\[
DP(0,0)=
\begin{pmatrix}
\sigma&0\\
0&\sigma
\end{pmatrix}.
\]
The eigenvalues of the return map are
therefore both equal to $e^{-aT}$.
The claimed stability conclusions follow.
\end{proof}

\section{Numerical explorations of the three-dimensional model}
\label{sec:numerics}

We now investigate the dynamics of the three-particle system
\eqref{eq:torus-angular-equations} numerically. Throughout this section,
we represent the phase space $\mathbb T^3$ by the identification cube
$[0,2\pi]^3$, with opposite faces identified.
We consider two choices of the effective matrices and use the
inverse-temperature parameter $\beta\geq0$ as a bifurcation parameter. 
Define the matrices
\begin{equation}\label{eq:numerics-main-matrices}
\begin{gathered}
A_1=\begin{pmatrix}0.20&-2.30\\2.30&0.20\end{pmatrix},
\qquad
A_2=\begin{pmatrix}-0.05&-2.30\\2.30&-0.05\end{pmatrix},\\[2pt]
B=\begin{pmatrix}-1.20&0.40\\-0.80&-2.40\end{pmatrix}.
\end{gathered}
\end{equation}
The first choice has $A=A_1$ and the second has $A=A_2$;
both use the same matrix $B$. For $A=A_1$, the synchronized
periodic orbit of Proposition~\ref{prop:torus-rotating-wave}
is attracting for every $\beta\geq0$ (see
Section~\ref{sec:attractingCycle}). Parameter scans initialized
near this orbit reveal no additional dynamics.

For $A=A_2$, the synchronized periodic orbit is transversely
repelling for every $\beta\geq0$, and a standard bifurcation
diagram reveals richer behavior. The diagram suggests torus-like
dynamics, periodic windows, and numerical evidence of chaos.
We consider these results in Section~\ref{sec:chaos1}.

For the same choice $A=A_2$, we also follow an attracting,
nonsynchronized periodic orbit, apparently born in a
saddle-node bifurcation near $\beta=0.40916$.
This periodic orbit has complex-conjugate transverse Floquet
multipliers, and continuation indicates a candidate
Neimark--Sacker bifurcation near $\beta=0.58556$.
We discuss this second scan in Section~\ref{sec:secondSweep}.

Our goal is not to provide a comprehensive bifurcation
analysis of this system. Rather, we exhibit several representative
dynamical phenomena, illustrating a wide variety of 
possible behaviors for only three interacting particles.
The three-particle computations use adaptive integration of the angular
equations, event detection for Poincar\'e crossings, and integration of
the variational equations. Periodic solutions are refined by shooting,
and their orbital stability is computed via Floquet multipliers.
Lyapunov exponents are estimated by repeated QR orthogonalization.
The coupled-system Lyapunov computations in Section~5 instead use
fixed-step fourth-order Runge--Kutta integration.
For general discussions of these numerical techniques, see
\cite{Doedel2007Lecture,Kuznetsov2004AppliedBifurcation,BenettinEtAl1980Numerical}.

Unless otherwise stated, our Poincar\'e section is
\[
\Sigma=\{\theta\in\mathbb T^3:\theta_1=0\pmod{2\pi}\},
\]
with coordinates $(\theta_2,\theta_3)$. This section is chosen for
convenience, rather than because it has any distinguished geometric
status. The bifurcation diagrams record the coordinate $\theta_3$
at section encounters, and the Poincar\'e maps distinguish crossings
with $\dot\theta_1<0$ from those with $\dot\theta_1>0$. As we shall see,
these two crossing directions can correspond to different portions
of the same geometric intersection curve when an invariant torus
is not everywhere transverse to the flow through $\Sigma$.

\subsection{An explicit attracting synchronized periodic orbit} \label{sec:attractingCycle}

As a first illustration, take $A=A_1$ and $B$ as in Equation
\eqref{eq:numerics-main-matrices}.
The synchronized circle $\Delta=\{\theta_1=\theta_2=\theta_3\}$
is invariant. By Proposition~\ref{prop:torus-rotating-wave}, it contains
the explicit periodic orbit
$\theta_1(t)=\theta_2(t)=\theta_3(t)=\varphi_0+2.30t\pmod{2\pi}$,
with period $T=2\pi/2.30$ and double nontrivial Floquet multiplier
$\exp(-0.20T)\approx0.57905$. In particular, the synchronized periodic
orbit is attracting. In fact, the synchronized circle
itself is this periodic orbit, and both its period and 
multipliers are independent of $\beta$ and $B$. Thus it is present
and normally attracting for every $\beta\geq0$, with $A=A_1$ fixed.

\begin{figure}[tbp]
    \centering
    \begin{subfigure}[t]{0.48\textwidth}
        \centering
        \includegraphics[width=\linewidth]{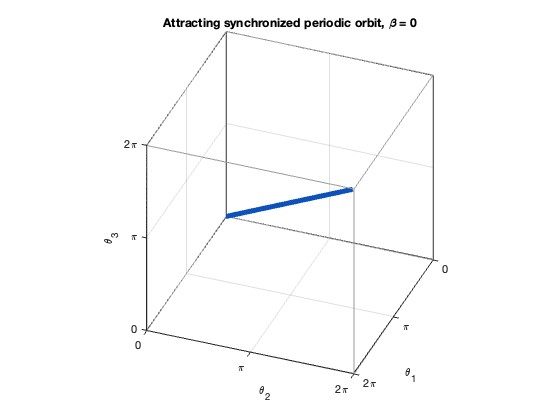}
        \caption{The attracting synchronized periodic orbit at $\beta=0$.}
        \label{fig:numerics-sync-orbit}
    \end{subfigure}
    \hfill
    \begin{subfigure}[t]{0.48\textwidth}
        \centering
        \includegraphics[width=\linewidth]{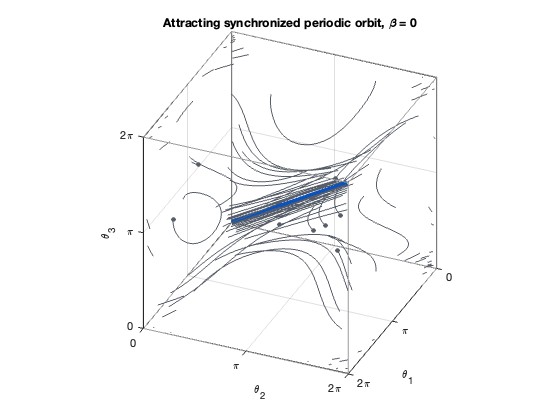}
        \caption{The same orbit together with selected forward trajectories, some of which approach it.}
        \label{fig:numerics-sync-transients}
    \end{subfigure}
    \caption{An explicitly known attracting periodic orbit in the
    synchronized circle (left) and the behavior of selected initial
    conditions (right).  The simulations use $A = A_1$ and the 
    matrix $B$ defined in Equation \eqref{eq:numerics-main-matrices}.
    The curves are drawn in the identification cube;
    trajectories are broken where they cross an identified face.}
    \label{fig:numerics-synchrony}
\end{figure}

Figure~\ref{fig:numerics-synchrony} illustrates the resulting
periodic motion and the approach of selected trajectories to the
synchrony circle. The numerical experiments also remind us that local
orbital stability does not determine the global basin of attraction.

\subsection{The main parameter scan: invariant tori and chaos} \label{sec:chaos1}

For the remaining experiments, we fix $A=A_2$ and the same
matrix $B$ from \eqref{eq:numerics-main-matrices}, varying $\beta\geq0$.
The matrices $A_1$ and $A_2$ have the same rotational part,
while their scalar parts are $0.20I_2$ and $-0.05I_2$,
respectively. The synchronized periodic orbit remains present for
$A=A_2$, but it is transversely repelling, since its two nontrivial
Floquet multipliers are
$\exp(0.05\,2\pi/2.30)\approx1.14636>1$.

Our first parameter scan begins at $\beta=0$ with the angular initial
condition
\[
 \theta(0)=(1.8647340967,\;3.8249352834,\;4.0754264806),
\]
and increases $\beta$ using the terminal state from one parameter
value to initialize the next. After discarding a transient, section
encounters are plotted against $\beta$. This procedure is well suited
to locating attracting regimes and sudden changes in the observed
asymptotic dynamics, but can also switch between coexisting attractors
or retain long transients.

\begin{figure}[tbp]
    \centering
    \includegraphics[width=0.80\textwidth]{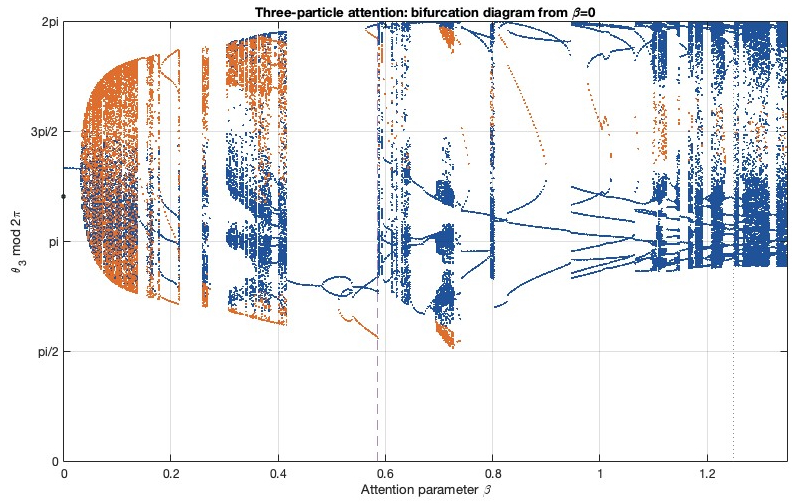}
    \caption{Bifurcation diagram for the main forward-integration scan
    with the matrices $A_2$ and $B$ from
    \eqref{eq:numerics-main-matrices}. Points record
    the $\theta_3$ coordinate at intersections with the section
    $\theta_1=0\pmod{2\pi}$ after a transient has been discarded.
    The diagram is an empirical representation of the attractors
    selected by this parameter sweep, rather than a continuation
    diagram of every periodic solution.}
    \label{fig:numerics-main-bifurcation}
\end{figure}

The diagram in Figure~\ref{fig:numerics-main-bifurcation} suggests a
rich progression of periodic, quasiperiodic, and more complicated
recurrent behavior. To examine individual regimes, we compute the
finite-time Lyapunov spectrum, a Poincar\'e section, and a phase-space
portrait at fixed values of $\beta$. The three panels in the following
figures show these diagnostics in that order. A smooth closed curve in
the section, together with two Lyapunov exponents approaching zero,
is consistent with an invariant two-torus. By contrast, a persistently
positive largest Lyapunov exponent is numerical evidence of sensitive
dependence on initial conditions.

\begin{figure}[tbp]
    \centering
    \begin{subfigure}[t]{0.32\textwidth}
        \centering
        \includegraphics[width=\linewidth]{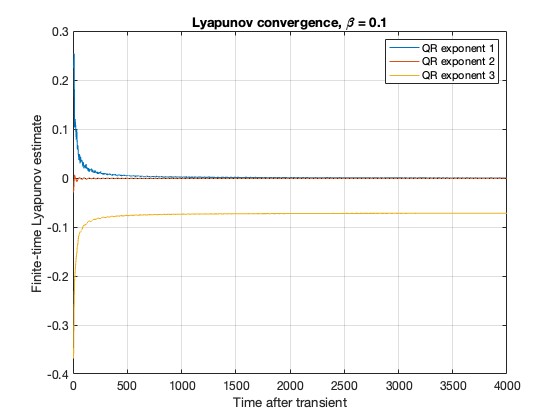}
        \caption{Finite-time Lyapunov exponents.}
        \label{fig:numerics-main-torus-lyap}
    \end{subfigure}
    \hfill
    \begin{subfigure}[t]{0.32\textwidth}
        \centering
        \includegraphics[width=\linewidth]{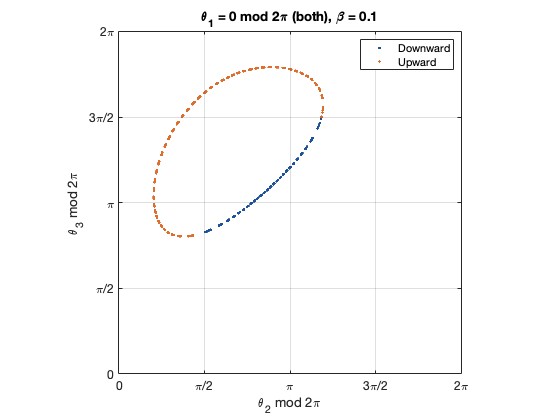}
        \caption{Poincar\'e intersection.}
        \label{fig:numerics-main-torus-section}
    \end{subfigure}
    \hfill
    \begin{subfigure}[t]{0.32\textwidth}
        \centering
        \includegraphics[width=\linewidth]{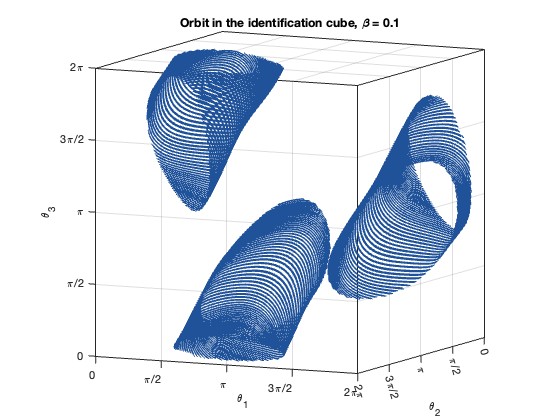}
        \caption{Orbit in the identification cube.}
        \label{fig:numerics-main-torus-phase}
    \end{subfigure}
    \caption{Dynamics of the main scan at $\beta=0.1$. The numerical
    diagnostics are consistent with an attracting invariant two-torus.
    In the middle panel, blue marks crossings with $\dot\theta_1<0$,
    while yellow/orange marks crossings with $\dot\theta_1>0$.
    The flow need not be everywhere transverse to the chosen section.}
    \label{fig:numerics-main-torus}
\end{figure}

The two colors in the Poincar\'e plots merit an explanation. The
blue points correspond to trajectories crossing $\Sigma$ with
$F_1(\theta)=\dot\theta_1<0$ and the yellow/orange points to crossings
with $F_1(\theta)>0$. When an invariant two-torus intersects $\Sigma$
along a closed curve, the vector field need not point to the same side
of $\Sigma$ along the entire curve. At transitions between these two
sets of intersection points, continuity forces $F_1=0$, so the vector
field is tangent to the chosen section. In particular, the complete
geometric intersection curve need not be a global transverse
Poincar\'e section for the flow. This geometric effect explains why
plotting only one crossing direction can reveal just an arc of the
closed intersection curve, even when the phase-space portrait clearly
suggests a smooth torus.

At $\beta=1.25$, the same matrices give a substantially different
picture. Numerical integrations from different initial conditions
yield similar Lyapunov spectra, with one positive exponent, one
exponent close to zero, and one negative exponent (approximately
$0.11$, $0$, and $-0.31$). Along with the nontrivial structure of the
Poincar\'e section, this provides evidence for chaotic dynamics,
although it does not constitute a proof of the existence of a
chaotic attractor.

\begin{figure}[tbp]
    \centering
    \begin{subfigure}[t]{0.32\textwidth}
        \centering
        \includegraphics[width=\linewidth]{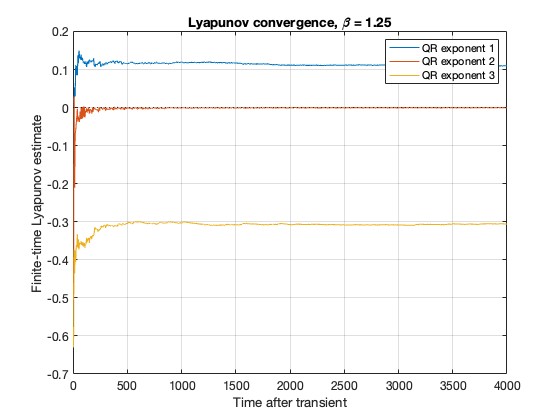}
        \caption{Finite-time Lyapunov exponents.}
        \label{fig:numerics-main-chaos-lyap}
    \end{subfigure}
    \hfill
    \begin{subfigure}[t]{0.32\textwidth}
        \centering
        \includegraphics[width=\linewidth]{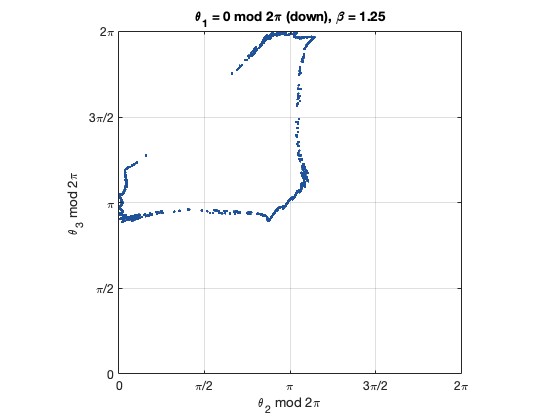}
        \caption{Poincar\'e intersection.}
        \label{fig:numerics-main-chaos-section}
    \end{subfigure}
    \hfill
    \begin{subfigure}[t]{0.32\textwidth}
        \centering
        \includegraphics[width=\linewidth]{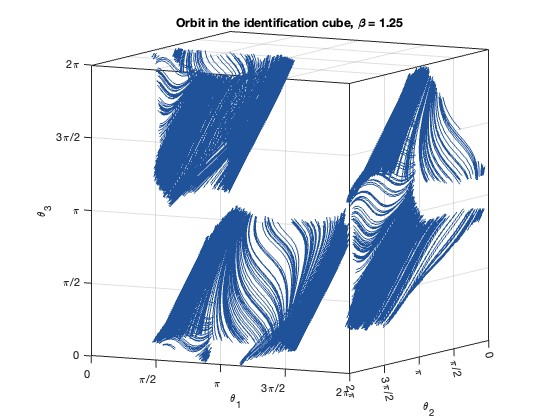}
        \caption{Orbit in the identification cube.}
        \label{fig:numerics-main-chaos-phase}
    \end{subfigure}
    \caption{Numerical evidence of chaotic dynamics along the main
    scan of the $A_2$ system
     at $\beta=1.25$. The positive largest Lyapunov exponent and
    the structure of the section distinguish this regime from the
    smooth torus illustrated in Figure~\ref{fig:numerics-main-torus}.}
    \label{fig:numerics-main-chaos}
\end{figure}

\subsection{A second parameter scan: a fold of periodic orbits and
complex Floquet multipliers} \label{sec:secondSweep}

A second set of experiments uses the same matrices $A=A_2$ and $B$
from~\eqref{eq:numerics-main-matrices}, but is initialized near a
nonsynchronized attracting periodic orbit at $\beta=0.55$.
A convenient section point on this periodic orbit is
\[
 \theta(0)\approx(0,\;3.4617585711,\;3.5395509760),
 \qquad T\approx13.79255247.
\]
Its two nontrivial Floquet multipliers are approximately
\[
 \mu_{1,2}=-0.290456\pm0.326888\,\mathrm{i},
\]
and hence have modulus smaller than one. Thus the orbit is
attracting, with a contracting two-dimensional normal
bundle on which the return dynamics has a rotational component.
This provides a qualitatively different example from the synchronized
orbit, whose nontrivial multipliers are real.

\begin{figure}[tbp]
    \centering
    \begin{subfigure}[t]{0.48\textwidth}
        \centering
        \includegraphics[width=\linewidth]{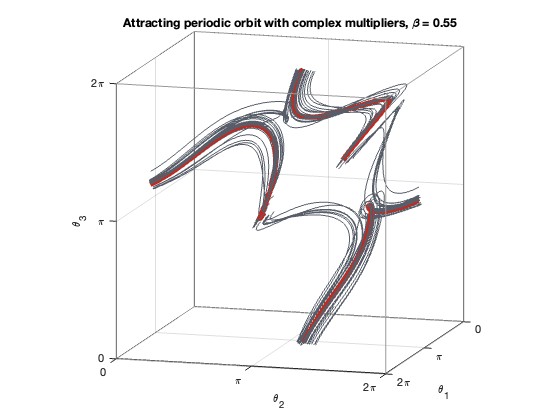}
        \caption{The nonsynchronized attracting periodic orbit at $\beta=0.55$.}
        \label{fig:numerics-complex-orbit}
    \end{subfigure}
    \hfill
    \begin{subfigure}[t]{0.48\textwidth}
        \centering
        \includegraphics[width=\linewidth]{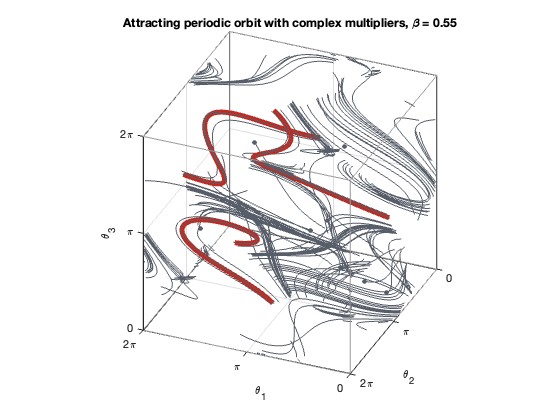}
        \caption{The same orbit with selected trajectories illustrating its basin of attraction.}
        \label{fig:numerics-complex-transients}
    \end{subfigure}
    \caption{An attracting nonsynchronized periodic orbit with
    complex-conjugate transverse Floquet multipliers for the $A_2$ system. 
    Several plotted
    trajectories do not approach this orbit on the displayed time scale,
    suggesting that its basin of attraction is not the whole phase
    space; the plots alone do not determine the basin boundary.}
    \label{fig:numerics-complex-periodic}
\end{figure}

Numerical shooting and continuation suggest that this periodic orbit
belongs to a family created in a saddle-node bifurcation of periodic
orbits at
\[
 \beta_{\mathrm{fold}}\approx0.40916194.
\]
At the candidate fold, one nontrivial Floquet multiplier is close to
$+1$, while the other is approximately $0.20636$. Immediately above
the fold, the calculations find both a transversely attracting
periodic orbit and a saddle periodic orbit. In particular, the
nonsynchronized family need not extend as a regular branch all the
way to $\beta=0$. Numerical continuation of this periodic-orbit family also
indicates a crossing of a complex-conjugate pair of multipliers
through the unit circle near
$\beta_{\mathrm{NS}}\approx0.58555751$, consistent with a
Neimark--Sacker bifurcation. Establishing the bifurcation type and the
local organization of invariant tori requires additional work.

The second bifurcation scan starts at the $\beta=0.55$ periodic-orbit
seed and varies $\beta$ in both directions. As in the main scan, the
points shown are obtained from forward integration, so their
appearance depends on which invariant set is reached from the chosen
initial conditions. The superimposed information from shooting
continuation, when displayed, must be interpreted separately from
the attracting-set samples.

\begin{figure}[tbp]
    \centering
    \includegraphics[width=1.0\textwidth]{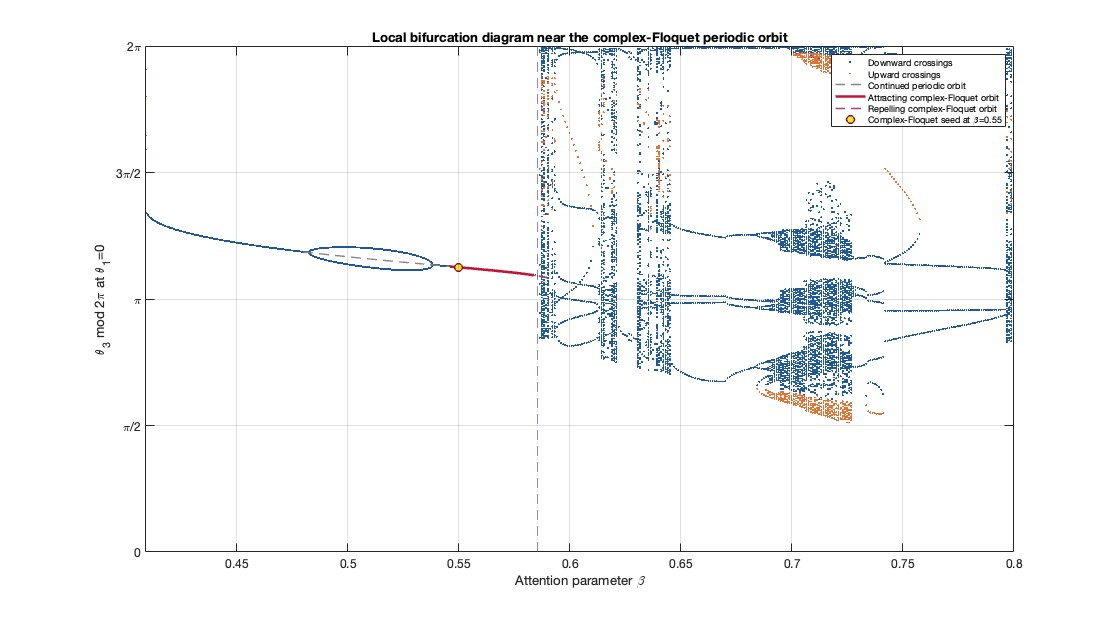}
    \caption{Bifurcation diagram from the second, complex-orbit-seeded
    parameter scan, using the same matrices as in
    Figure~\ref{fig:numerics-main-bifurcation}. Forward-orbit
    Poincar\'e points reveal the asymptotic dynamics selected near
    the nonsynchronized periodic family. Any continued periodic-orbit
    branch shown in the plot is distinct from the point cloud generated
    by time integration.}
    \label{fig:numerics-complex-bifurcation}
\end{figure}

The second scan illustrates a further subtlety of the dynamics:
nearby parameter values may support different recurrent sets.
For example, continuation indicates that the nonsynchronized periodic
orbit is still attracting at $\beta=0.58$, while a separate
forward-integration experiment produces a torus-like portrait at
$\beta=0.59$. These observations are consistent with a change in
stability near the candidate Neimark--Sacker value, but do not by
themselves establish the bifurcation mechanism or coexistence of
attractors at a fixed parameter value.

\begin{figure}[tbp]
    \centering
    \begin{subfigure}[t]{0.32\textwidth}
        \centering
        \includegraphics[width=\linewidth]{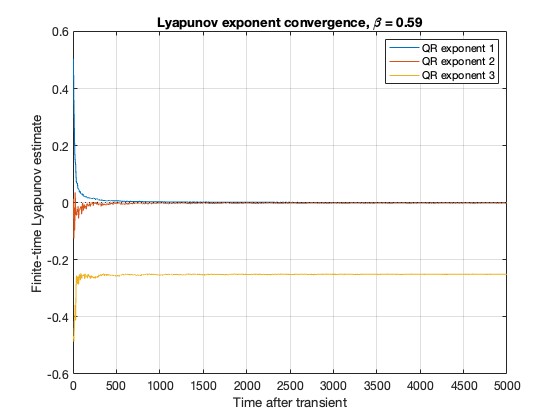}
        \caption{Finite-time Lyapunov exponents.}
        \label{fig:numerics-complex-torus-lyap}
    \end{subfigure}
    \hfill
    \begin{subfigure}[t]{0.32\textwidth}
        \centering
        \includegraphics[width=\linewidth]{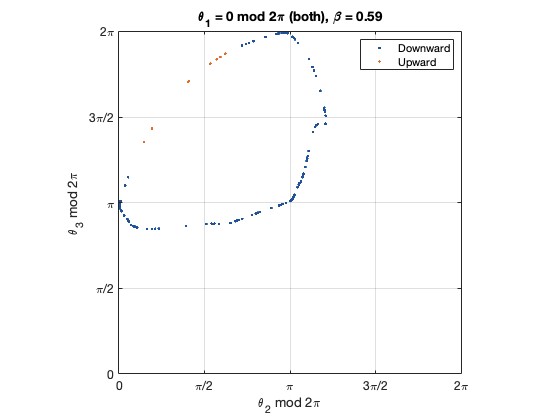}
        \caption{Poincar\'e intersection.}
        \label{fig:numerics-complex-torus-section}
    \end{subfigure}
    \hfill
    \begin{subfigure}[t]{0.32\textwidth}
        \centering
        \includegraphics[width=\linewidth]{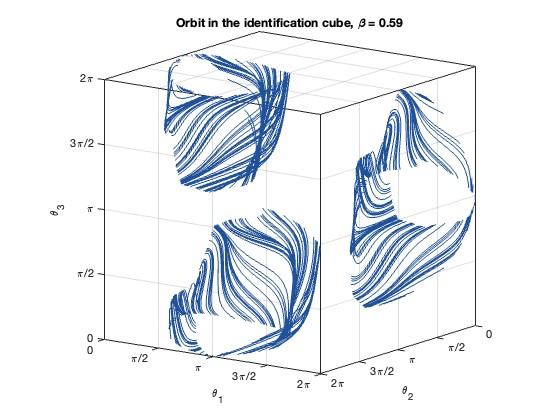}
        \caption{Orbit in the identification cube.}
        \label{fig:numerics-complex-torus-phase}
    \end{subfigure}
    \caption{A torus-like recurrent set observed in the second scan
    at $\beta=0.59$. The numerical portrait does not by itself
    establish the existence of an invariant torus or its
    bifurcation mechanism.}
    \label{fig:numerics-complex-torus}
\end{figure}

\begin{figure}[tbp]
    \centering
    \begin{subfigure}[t]{0.32\textwidth}
        \centering
        \includegraphics[width=\linewidth]{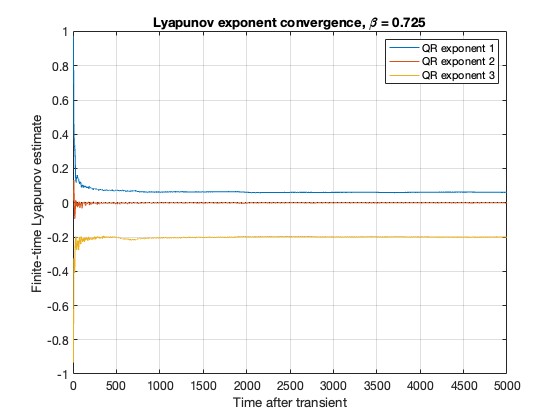}
        \caption{Finite-time Lyapunov exponents.}
        \label{fig:numerics-complex-higher-lyap}
    \end{subfigure}
    \hfill
    \begin{subfigure}[t]{0.32\textwidth}
        \centering
        \includegraphics[width=\linewidth]{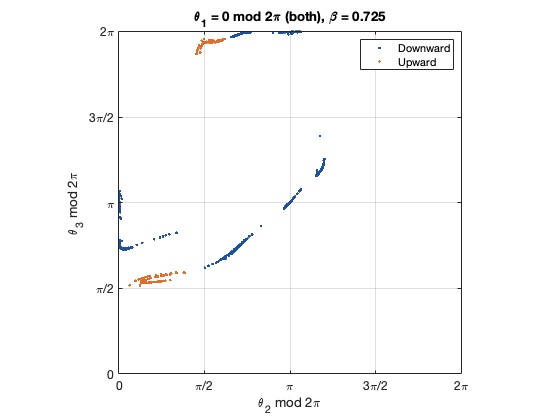}
        \caption{Poincar\'e intersection.}
        \label{fig:numerics-complex-higher-section}
    \end{subfigure}
    \hfill
    \begin{subfigure}[t]{0.32\textwidth}
        \centering
        \includegraphics[width=\linewidth]{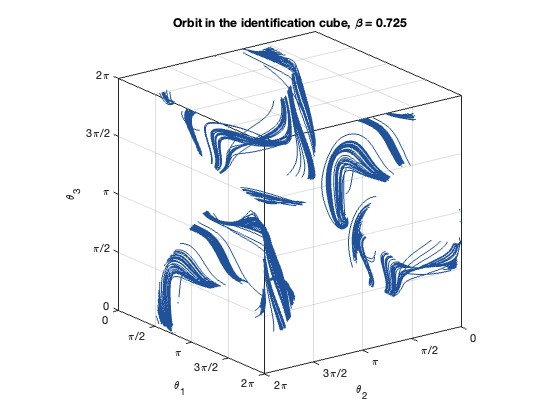}
        \caption{Orbit in the identification cube.}
        \label{fig:numerics-complex-higher-phase}
    \end{subfigure}
    \caption{Numerical evidence of chaotic dynamics in the second
    parameter scan at $\beta=0.725$. The finite-time Lyapunov spectrum
    (left), Poincar\'e section (middle), and phase-space portrait
    (right) provide complementary diagnostics of the chaotic regime.}
    \label{fig:numerics-complex-higher}
\end{figure}

At $\beta=0.725$, the second scan exhibits a markedly more
complicated recurrent regime. The finite-time Lyapunov computation
for the trajectory illustrated in Figure~\ref{fig:numerics-complex-higher}
indicates a positive largest exponent. Together with the structure
of the Poincar\'e section and the phase-space portrait, this provides
numerical evidence for chaotic dynamics. This example is especially
interesting when viewed alongside the attracting periodic orbit at
$\beta=0.55$ and the torus-like dynamics at $\beta=0.59$.
The computations suggest that a single family of attention models
can exhibit several distinct kinds of recurrent behavior, although
they do not constitute a rigorous proof of a chaotic invariant set.

Taken together, these computations illustrate that even the
three-particle projected-attention model can exhibit a variety of
recurrent phenomena. The apparent fold of nonsynchronized periodic
orbits, the emergence of complex Floquet multipliers, and the possible
interactions among periodic orbits, invariant tori, and chaotic invariant
sets provide concrete targets for more systematic numerical continuation
and, ultimately, computer-assisted validation.

\section{Higher-dimensional transformers as coupled systems}
\label{sec:bigTransformers}

The three-particle system studied in Section~\ref{sec:torus} 
will now serve as a component for assembling larger systems
with interesting dynamics. Consider $M\geq2$ groups of three
particles, with index sets $G_1,\ldots,G_M$ partitioning
$\{1,\ldots,3M\}$.  Suppose initially that all particles lie on the
same ellipse $\mathcal E_W$ and that, within each group, they evolve according
to the autonomous three-particle system derived in
Section~\ref{sec:torus}.  Writing $X^{(m)}=(x_i)_{i\in G_m}$,
the uncoupled equations are
\begin{equation}\label{eq:big-product}
 \dot X^{(m)}=F_3(X^{(m)}),\qquad m=1,\ldots,M.
\end{equation}
The phase space is $\mathcal E_W^{3M}\cong\mathbb T^{3M}$, and the
flow of \eqref{eq:big-product} is the direct product of the constituent
three-particle flows.  In particular, recurrent invariant sets of the
individual models become invariant sets of the larger system
by taking products.

We now introduce a coupling parameter $\varepsilon\geq0$ and a
family of projected attention vector fields of the form
\begin{equation}\label{eq:big-perturbation}
 \dot X=F^\varepsilon(X),\qquad F^0(X)
 =\bigl(F_3(X^{(1)}),\ldots,F_3(X^{(M)})\bigr),
\end{equation}
with $F^\varepsilon\to F^0$ in $C^r$ for each fixed $M$
and finite $r$, as $\varepsilon\to0$. 
In the next subsections we discuss several ways to 
assemble a self-attention network of this form, using 
different techniques to retain subsystem information.

Note that we must build subsystem information into the network 
if we want to treat self-attention transformers as coupled oscillators.  
This is because ordinary, unmasked attention, as in 
Equation~\eqref{eq:intro-particle-equations},
is equivariant under permutations of the $3M$ particles.  A
nontrivial direct product of independent three-particle systems, by
contrast, distinguishes the prescribed subgroups, and
permutations that mix particles from different groups
actually change the dynamics.  We will see that masks and fixed group labels provide 
the structural information we need.

There are numerous precedents for structured interactions in machine learning.
Graph attention networks restrict interactions by specifying a prescribed graph
\cite{VelickovicEtAl2018}; heterogeneous graph transformers permit
node- and edge-type-dependent interactions \cite{HuEtAl2020}; and
head-specific masks, as discussed below, are studied in  
\cite{WangEtAl2020RoleMasks}.  
We now describe three realizations, all of which recover
Equation~\eqref{eq:big-product} in the uncoupled limit.
Two use masks, and one uses group labels. For positive
coupling, the first two realizations define exactly the same
coupled vector field, while the third uses a different
normalization architecture.

\paragraph{Multiplicative attention masks.}
We use a common convention for attention masks throughout
this section. Given a collection of scores $s_{ij}(X)$,
a nonnegative mask $m=(m_{ij})$ defines the modified
attention coefficients
\[
\alpha_{ij}(X;m)
=
\frac{m_{ij}e^{s_{ij}(X)}}
{\displaystyle\sum_{k=1}^{N}
m_{ik}e^{s_{ik}(X)}}.
\]
We assume that each row of $m$ contains at least one
positive entry, so that the denominator is nonzero.
A hard mask has entries in $\{0,1\}$, while more general
nonnegative masks allow selected interactions to be
weakened or strengthened.

Consider $N=3M$ particles in $\mathbb{R}^2$,
partitioned into $M\geq 2$ disjoint groups, each
containing three particles. Write
\[
\{1,\ldots,3M\}=G_1\sqcup\cdots\sqcup G_M.
\]
We introduce the complementary hard masks
\[
b_{ij}^{\rm int}
=
\begin{cases}
1,&i,j\in G_m\text{ for some }m,\\
0,&\text{otherwise},
\end{cases}
\qquad
b_{ij}^{\rm ext}=1-b_{ij}^{\rm int}.
\]
Thus $b^{\rm int}$ permits only interactions between
particles in the same group, while $b^{\rm ext}$
permits only interactions between particles belonging
to distinct groups.

We consider two ways of using these masks to construct
weakly coupled attention systems. In the first, the
internal and external masks are combined before applying
a single softmax normalization over all particles.
In the second, the two masks are applied to separate
attention heads, each with its own normalization, and
the resulting outputs are combined afterward.

These constructions agree in the uncoupled limit but
differ for positive coupling. In particular, in the
single-head construction the relative contribution of
external interactions depends on the number of particles
outside each group. In the two-head construction, the
external head produces a separately normalized weighted
average, whose magnitude remains bounded independently
of the number of groups. This distinction becomes
important when considering systems with many particles.

\subsection{Block-masked single-head attention}
\label{subsec:big-blockmask}

We first consider a single-head attention system in which
the constituent subsystems are identical. To this end,
fix common matrices $Q,K,V$ and an inverse temperature
$\beta\geq0$, and write
\[
s_{ij}(X)=\beta\langle Qx_i,Kx_j\rangle.
\]
Using the complementary hard masks $b^{\rm int}$ and
$b^{\rm ext}$ defined above, introduce the
parameter-dependent mask
\[
m_{ij}(\varepsilon)
=b_{ij}^{\rm int}+\varepsilon b_{ij}^{\rm ext},
\qquad 0\leq\varepsilon\leq1.
\]
Thus interactions within each group have weight one,
while interactions between distinct groups have
weight $\varepsilon$. The corresponding masked
attention coefficients are
\begin{equation}\label{eq:big-softmask}
\alpha_{ij}^\varepsilon(X)
=
\frac{m_{ij}(\varepsilon)e^{s_{ij}(X)}}
{\displaystyle\sum_{k=1}^{3M}
m_{ik}(\varepsilon)e^{s_{ik}(X)}},
\end{equation}
and the resulting projected equations are
\begin{equation}\label{eq:big-masked-field}
\dot x_i
=
\Pi^W_{x_i}\left(
\sum_{j=1}^{3M}\alpha_{ij}^\varepsilon(X)Vx_j
\right).
\end{equation}
Observe that this construction uses a single
softmax normalization over all particles,
rather than separate normalizations for
within-group and between-group interactions.

At $\varepsilon=0$, the mask is precisely
$b^{\rm int}$, so all interactions between
distinct groups vanish and each attention
row is normalized over its own group.
Consequently, \eqref{eq:big-masked-field}
reduces exactly to the uncoupled product
system \eqref{eq:big-product}.
At $\varepsilon=1$, all entries of the mask
are equal to one, and the model becomes
ordinary unmasked single-head attention
involving all $3M$ particles. For
$0<\varepsilon<1$, cross-group interactions
are suppressed but remain nonzero.

The dependence on $\varepsilon$ is smooth,
including at $\varepsilon=0$. Indeed, for
$i\in G_m$, define
\[
Z_i^{\rm in}=\sum_{j\in G_m}e^{s_{ij}},
\qquad
Z_i^{\rm out}=\sum_{j\notin G_m}e^{s_{ij}},
\]
and
\[
Y_i^{\rm in}=\sum_{j\in G_m}e^{s_{ij}}Vx_j,
\qquad
Y_i^{\rm out}=\sum_{j\notin G_m}e^{s_{ij}}Vx_j.
\]
The $i$th component of the vector field
can then be written as
\begin{equation}\label{eq:big-softmask-expansion}
\begin{aligned}
F_i^\varepsilon(X)
&=\Pi^W_{x_i}
\frac{Y_i^{\rm in}+\varepsilon Y_i^{\rm out}}
{Z_i^{\rm in}+\varepsilon Z_i^{\rm out}}\\
&=F_i^0(X)
+\varepsilon\Pi^W_{x_i}\left(
\frac{Y_i^{\rm out}}{Z_i^{\rm in}}
-\frac{Z_i^{\rm out}Y_i^{\rm in}}
{(Z_i^{\rm in})^2}
\right)
+O(\varepsilon^2).
\end{aligned}
\end{equation}
Since $Z_i^{\rm in}$ is strictly positive,
and the configuration space is compact,
the remainder is uniform, together with
any fixed number of derivatives in $X$.
Thus, for each fixed number of groups and
each finite $r$,
\[
\|F^\varepsilon-F^0\|_{C^r}=O(\varepsilon).
\]
In particular, the block-masked system
is a smooth perturbation of the uncoupled
product system, so that standard persistence
results for hyperbolic invariant sets and
normally hyperbolic invariant manifolds
are applicable whenever their respective
hypotheses are satisfied.

For systems with many groups, the effective
strength of cross-group attention depends
not only on $\varepsilon$ but also on the
number of particles outside each group.
In particular, if the attention scores are
uniformly bounded independently of $M$,
then $Z_i^{\rm out}/Z_i^{\rm in}$ is of order
$M-1$, since the external sum has $3(M-1)$
terms and the internal sum has three. Thus
the effective relative contribution of cross-group
attention is controlled by $\varepsilon
Z_i^{\rm out}/Z_i^{\rm in}$. Consequently,
a fixed small value of $\varepsilon$ need not
represent equally weak coupling as the
number of groups increases. A natural choice
is to take $\varepsilon=O(M^{-1})$, with a
sufficiently small proportionality constant,
or to restrict cross-group attention to a
sparse network. This scaling consideration
also applies to the equivalent group-label
construction below, but not to the separately
normalized two-head construction of
Section~\ref{subsec:masked-coupled-systems}.

The next subsection gives an equivalent
implementation using fixed group labels.
Section~\ref{subsec:masked-coupled-systems}
then introduces the alternative two-head construction.

\subsection{Group labels and attention biases}
\label{subsec:big-labels}

The preceding mask can also be implemented
using ordinary softmax attention applied to \emph{augmented}
particle representations.  That is, one assigns to each group $G_m$, 
$1 \leq m \leq M$  a fixed label
$e_m\in\mathbb R^M$, the $m$th coordinate vector.   When we 
append a label to each particle, we obtain $\widehat x_i=(x_i,e_m)$ 
for $i\in G_m$.  Define augmented query and key maps with 
bilinear compatibility score defined by 
\begin{equation}\label{eq:big-label-score}
 \widehat s_{ij}
 =\beta\langle Qx_i,Kx_j\rangle
 +\gamma\langle e_{g(i)},e_{g(j)}\rangle
 =s_{ij}+\gamma\,\mathbf1_{\{g(i)=g(j)\}},
\end{equation}
where $g(i)$ denotes the group containing particle $i$.  Such a
score is obtained by a block-diagonal choice of
$\widehat Q^\top\widehat K$ on the augmented space.  The values
remain $Vx_j$, and the labels are fixed rather than evolving
coordinates.  Factoring $e^\gamma$ from each softmax denominator
shows that the resulting attention coefficients coincide exactly
with \eqref{eq:big-softmask}, with
\begin{equation}\label{eq:big-label-epsilon}
 \varepsilon=e^{-\gamma}.
\end{equation}
Hence, arbitrarily weak positive coupling is implemented using finite
query and key matrices, without an attention mask.  Exact
decoupling corresponds to the limit $\gamma\to+\infty$; it cannot
be attained by a finite value of this particular bias parameter.

This formulation makes the group identity part of the token
representation and is closely related to type-dependent attention
and additive structural biases: score modifications supported
by practical attention software \cite{PyTorch2024FlexAttention}.  
It is not a construction using only unaugmented positions on the ellipse: 
the additional fixed labels distinguish the subsystems.

\subsection{Two-head masked attention as a coupled system}
\label{subsec:masked-coupled-systems}

The second normalization architecture, perhaps the most useful
for the purposes of the present work, uses two masked attention heads.
The notion of multi-head attention originates in
\cite{VaswaniEtAl2017}, and we recall the definition for
two heads before specializing to a masked partition of
particles into interacting groups.

\paragraph{Two-head self-attention.}
Let $X=(x_1,\ldots,x_N)$ denote a configuration of
$N$ particles in a manifold $\mathcal M\subset\mathbb R^d$.
In a two-head self-attention mechanism, each head
$h\in\{1,2\}$ has its own query, key, and value matrices
$Q_h$, $K_h$, and $V_h$, together with an inverse
temperature parameter $\beta_h\geq0$.

The attention weights associated with head $h$ are
\[
\rho_{ij}^{(h)}(X)
=
\frac{
\exp\bigl(\beta_h\langle Q_hx_i,K_hx_j\rangle\bigr)
}{
\displaystyle\sum_{k=1}^N
\exp\bigl(\beta_h\langle Q_hx_i,K_hx_k\rangle\bigr)
},
\qquad 1\leq i,j\leq N.
\]
Each head produces an attention-weighted average
\[
z_i^{(h)}(X)
=
\sum_{j=1}^N\rho_{ij}^{(h)}(X)V_hx_j.
\]
The two head outputs are combined through output
matrices $O_1$ and $O_2$, giving
\[
\mathcal A_i(X)
=
O_1z_i^{(1)}(X)+O_2z_i^{(2)}(X).
\]
This is equivalent to concatenating the two head
outputs and applying a single linear output map.
The corresponding continuous-time projected
attention system is
\[
\dot x_i=\Pi^W_{x_i}\mathcal A_i(X),
\qquad 1\leq i\leq N,
\]
where $\Pi^W_{x_i}$ denotes the tangent projection
associated with the chosen particle manifold.

\paragraph{Masked two-head attention.}
We now introduce a mask restricting which particles 
contribute to each attention head.  Using the
multiplicative-mask convention introduced at the
beginning of this section, we say that for each $h\in\{1,2\}$,
$m^{(h)}$ is a nonnegative mask matrix, with
at least one positive entry in every row. The
masked attention weights are defined by
\[
\rho_{ij}^{(h)}(X;m^{(h)})
=
\frac{
m_{ij}^{(h)}
\exp\bigl(\beta_h\langle Q_hx_i,K_hx_j\rangle\bigr)
}{
\displaystyle\sum_{k=1}^N
m_{ik}^{(h)}
\exp\bigl(\beta_h\langle Q_hx_i,K_hx_k\rangle\bigr)
}.
\]
In particular, a hard mask has entries in
$\{0,1\}$, with an entry $m_{ij}^{(h)}=0$
excluding particle $j$ completely from the
attention of particle $i$ in head $h$.

Equivalently, if $S_h(i)$ denotes the set
of particles visible to particle $i$ in
head $h$, then for a hard mask
\[
\rho_{ij}^{(h)}(X;m^{(h)})
=
\begin{cases}
\displaystyle
\frac{
e^{\beta_h\langle Q_hx_i,K_hx_j\rangle}
}{
\sum_{k\in S_h(i)}
e^{\beta_h\langle Q_hx_i,K_hx_k\rangle}
},
&j\in S_h(i),\\[8pt]
0,&j\notin S_h(i).
\end{cases}
\]
Note that the attention weights are normalized
separately within each head, so that
\[
\sum_{j=1}^N\rho_{ij}^{(h)}(X;m^{(h)})=1,
\qquad h\in\{1,2\}.
\]
Thus each head computes a weighted average
over its permitted particles, and the two
averages are combined only afterward.

\paragraph{Complementary group masks.}
Suppose that the particle indices are partitioned
into $M\geq2$ nonempty, disjoint groups,
\[
\{1,\ldots,N\}
=
G_1\sqcup\cdots\sqcup G_M.
\]
In our numerical experiments, each group contains
three particles, so that $N=3M$.

We introduce two attention heads, referred to
as the internal and external heads. The internal
head permits interactions only between particles
belonging to the same group, while the external
head permits interactions only between particles
belonging to different groups.

More precisely, we use the complementary hard
masks $b^{\rm int}$ and $b^{\rm ext}$ defined
at the beginning of this section, setting
\[
m^{(1)}=b^{\rm int},
\qquad
m^{(2)}=b^{\rm ext}.
\]
Thus every pair of particles is visible to
exactly one of the two heads. In particular,
self-interactions belong to the internal head.

Let $i\in G_m$. Applying the internal mask
to the first attention head gives
\[
\rho_{ij}^{\rm int}(X)
=
\frac{
e^{\beta_{\rm int}
\langle Q_{\rm int}x_i,K_{\rm int}x_j\rangle}
}{
\displaystyle\sum_{k\in G_m}
e^{\beta_{\rm int}
\langle Q_{\rm int}x_i,K_{\rm int}x_k\rangle}
},
\qquad j\in G_m,
\]
and $\rho_{ij}^{\rm int}=0$ for $j\notin G_m$.

Similarly, the external mask gives
\[
\rho_{ij}^{\rm ext}(X)
=
\frac{
e^{\beta_{\rm ext}
\langle Q_{\rm ext}x_i,K_{\rm ext}x_j\rangle}
}{
\displaystyle\sum_{k\notin G_m}
e^{\beta_{\rm ext}
\langle Q_{\rm ext}x_i,K_{\rm ext}x_k\rangle}
},
\qquad j\notin G_m,
\]
and $\rho_{ij}^{\rm ext}=0$ for $j\in G_m$.

The outputs of the two heads are therefore
\[
z_i^{\rm int}(X)
=
\sum_{j\in G_m}\rho_{ij}^{\rm int}(X)
V_{\rm int}x_j,
\qquad
z_i^{\rm ext}(X)
=
\sum_{j\notin G_m}\rho_{ij}^{\rm ext}(X)
V_{\rm ext}x_j.
\]

We now choose the output matrices to be
\[
O_{\rm int}=I_d,
\qquad
O_{\rm ext}=\varepsilon I_d,
\]
where $\varepsilon\geq0$ controls the strength
of the inter-group coupling. Substituting
these choices into the two-head projected
attention system gives
\begin{equation}\label{eq:big-twoheads}
\dot x_i
=
\Pi^W_{x_i}\left(
\sum_{j\in G_m}\rho_{ij}^{\rm int}(X)
V_{\rm int}x_j
+
\varepsilon
\sum_{j\notin G_m}\rho_{ij}^{\rm ext}(X)
V_{\rm ext}x_j
\right),
\qquad i\in G_m.
\end{equation}
This is precisely the weakly coupled particle
system considered below.

Writing $X^{(m)}=(x_i)_{i\in G_m}$, and using
the linearity of the tangent projection in its
vector argument, we may express the system as
\[
\dot X^{(m)}
=
F^{(m)}(X^{(m)})
+
\varepsilon H_m(X),
\qquad 1\leq m\leq M,
\]
where $F^{(m)}$ is the internal attention
vector field on the $m$th group and $H_m$
collects the external interactions.

In particular, when $\varepsilon=0$, the
external head makes no contribution and the
system reduces exactly to a product of
uncoupled subsystems,
\[
\dot X^{(m)}=F^{(m)}(X^{(m)}),
\qquad 1\leq m\leq M.
\]
Since the attention weights of the two heads
are normalized independently, neither
$F^{(m)}$ nor $H_m$ depends on $\varepsilon$.
Consequently, the full vector field depends
affinely on the coupling parameter, and
\[
\|F^\varepsilon-F^0\|_{C^r}
=O(\varepsilon)
\]
for each fixed number of groups and each
finite $r$.

For the constructions considered here, all groups
have identical internal attention parameters,
and the external attention parameters are also
independent of the interacting groups. Hence
the uncoupled subsystems are identical, although
their states and trajectories need not be.
These uniformity assumptions are convenient
for studying the effect of coupling but are
not intrinsic to the masked attention construction.

Unlike the block-masked single-head construction
of Section~\ref{subsec:big-blockmask}, the present
system uses separate softmax normalizations for
the internal and external heads. In particular,
the coupling parameter scales the output of
the external head without altering the internal
attention weights. Thus the two constructions
have the same uncoupled limit, but define
different coupled systems for $\varepsilon>0$.

\begin{remark}[Heterogeneous coupled subsystems] \label{rem:hetCoupledSys}
The assumption that all groups have identical internal
attention parameters, and that all inter-group interactions
are governed by the same coupling parameters, is made
for simplicity. More generally, one could allow the
matrices $Q_{\rm int}^{(m)}$, $K_{\rm int}^{(m)}$,
and $V_{\rm int}^{(m)}$, as well as the inverse
temperature $\beta_{\rm int}^{(m)}$, to depend
on the group index $m$. Similarly, the external
attention parameters could depend pairwise on
the interacting groups.
This generalization permits the coupling of
subsystems exhibiting qualitatively different dynamics,
including periodic, quasiperiodic, and chaotic behavior.
\end{remark}

\subsection{Perturbative persistence}

All three implementations give smooth perturbations of
a direct-product flow for each fixed number of groups.
Consequently, the classical theory of
normally hyperbolic invariant manifolds provides a natural
framework for understanding the persistence of certain
dynamical structures under weak coupling. This theory applies
not only to hyperbolic periodic orbits, but also to normally
hyperbolic invariant tori and, more generally, to compact
normally hyperbolic invariant manifolds of arbitrary dimension.
Such manifolds may be attracting, repelling, or of saddle type.
We record the following consequence.

\begin{proposition}
\label{prop:big-periodic-product}
Suppose that the $m$th three-particle subsystem admits a
compact, smooth, normally hyperbolic invariant manifold
$\mathcal N_m$, for $1\leq m\leq M$. Assume moreover that
the normal contraction and expansion rates of the individual
manifolds dominate the tangential growth rates of all
factors, so that the product
\[
\mathcal N_0
=
\mathcal N_1\times\cdots\times\mathcal N_M
\]
is a normally hyperbolic invariant manifold for the uncoupled
system \eqref{eq:big-product}. Then, under any of the above
coupling implementations, for all sufficiently small
$\varepsilon\geq 0$, there exists a nearby normally hyperbolic
invariant manifold $\mathcal N_\varepsilon$, diffeomorphic
to $\mathcal N_0$. The dimensions of its stable and unstable
normal bundles are preserved under the perturbation.
In particular, if $\mathcal N_0$ is normally attracting,
then so is $\mathcal N_\varepsilon$.
\end{proposition}

\begin{proof}
At $\varepsilon=0$, the derivative of the product flow is
block diagonal, and the invariant splitting along
$\mathcal N_0$ is obtained by taking the direct sums of
the corresponding stable, tangent, and unstable bundles
of the individual factors. The assumed domination of
normal rates over tangential rates ensures that this
splitting is normally hyperbolic. Since the coupling
constructions converge to the product vector field in
$C^1$ as $\varepsilon\to 0$, the conclusion follows
from the persistence theorem for normally hyperbolic
invariant manifolds
\cite{Fenichel1972,HirschPughShub1977}.
\end{proof}

An important special case occurs when each $\mathcal N_m$
is a hyperbolic periodic orbit. In this case the tangent
dynamics have zero exponential growth rates, while all
transverse directions are exponentially contracting or
expanding. The product is therefore a normally hyperbolic
invariant $M$-torus, and the proposition guarantees its
persistence under weak coupling. Similar conclusions hold
for products of normally hyperbolic invariant tori, provided
the required spectral gaps are satisfied. The dynamics
on the perturbed manifold need not be conjugate to the
original dynamics. For example, weak coupling between
periodic oscillators may produce phase locking rather than
quasiperiodic motion. The proposition does not specify a numerical
threshold for $\varepsilon$; in particular, it does not by itself
justify persistence at $\varepsilon=0.1$, as explored in Section~5.

The theory of persistence extends beyond smooth invariant
manifolds. Classical structural stability results for compact,
locally maximal uniformly hyperbolic invariant sets imply
that such sets persist under sufficiently small perturbations,
with their dynamics preserved up to topological conjugacy
(or orbit equivalence for flows)
\cite{Shub1987,KatokHasselblatt1995}. 
Note that these results do not
establish the persistence of the chaotic
dynamics observed in this paper. Numerical evidence of
positive Lyapunov exponents does not establish uniform
hyperbolicity, and a product of hyperbolic invariant sets
for independent flows generally has additional neutral
directions associated with the independent time
translations of the factors. Such products need not
themselves be uniformly hyperbolic.

We view these observations as motivation for further work
rather than as limitations of the dynamical systems
approach. The persistence theory suggests that a wide
variety of recurrent invariant sets may survive the
introduction of weak attention coupling, but determining
which structures actually arise in transformer models
remains an open problem. Of particular interest is the
possibility that metastable clustering may be organized
by the stable manifolds of saddle-type invariant sets
more complicated than equilibria, including periodic
orbits, invariant tori, and chaotic invariant sets.
Their slow stable and unstable directions, together
with the geometry of their invariant bundles, may
provide a framework for understanding long transient
behavior even in the absence of gradient structure.

The examples developed here suggest that the dynamical
possibilities of transformer particle systems are
considerably richer than those described by convergence
to consensus alone. A systematic investigation of
hyperbolicity, normal hyperbolicity, bifurcations,
and the persistence of invariant structures under
coupling is therefore a promising direction for
future research.

\subsection{Structured attention in practical architectures and software}
\label{sec:structureApplications}
The structured interactions considered above are motivated in part
by mechanisms already used in practical transformer architectures.
For example, Mistral AI's Mistral 7B language model uses sliding-window
attention to restrict the preceding tokens attended to by each layer
\cite{JiangEtAl2023Mistral,MistralAI2023Announcement}, while Google's
Gemma 2 alternates local and global attention layers
\cite{JiKumar2024Gemma2}. The Swin Transformer uses local and shifted
windows in computer vision, including object detection and image
segmentation \cite{LiuEtAl2021Swin}. These windowed architectures
illustrate the use of restricted interactions, although their masks
do not directly implement the weakly coupled blocks studied here.

Tencent's CCFormer provides a particularly relevant industrial
example: its recommendation architecture combines different
classes of information through feature-field-separated cross-attention
with relation-specific attention pathways and output fusion.
Its authors report deployment in production video-recommendation
and advertising-ranking systems \cite{WangEtAl2026CCFormer}.
The architectural principle of separately normalized interactions
is close to that of our two-head construction, though the detailed
models are different. For comparison, Alibaba's deployed HetMatch
system uses type-dependent interactions for sponsored-search
keyword matching, but it is a heterogeneous graph neural network
rather than a Transformer attention model
\cite{LiuEtAl2021HetMatch}.

The use of group labels and attention-score biases also has
precedents in practical architectures. BERT employs segment
embeddings \cite{DevlinEtAl2019BERT}, and relative-position
and linear attention biases are used in other transformer models
\cite{ShawEtAl2018Relative,PressEtAl2022ALiBi}. These mechanisms
are not identical to the group-label implementation above, but
show that group-identifying embeddings and modified attention
scores are familiar architectural ingredients. On the software side, PyTorch's
FlexAttention provides programmable attention-score modifications
and block masks, including head-dependent masks
\cite{PyTorch2024FlexAttention}. NVIDIA's cuDNN attention library
also supports structured masking and additive score biases
\cite{NVIDIA2026CuDNNAttention}. Together these examples
illustrate that the structural ingredients used in our coupling
constructions are compatible with established attention architectures
and widely used software, although none should be interpreted as
an implementation of the particular dynamical systems studied here.

\section{Dynamics of some many particle systems} \label{sec:simulations}

\subsection{Coupling between two three-particle systems}
\label{sec:sixDim}

\begin{figure}[tbp]
    \centering
    \includegraphics[width=1.0\textwidth]{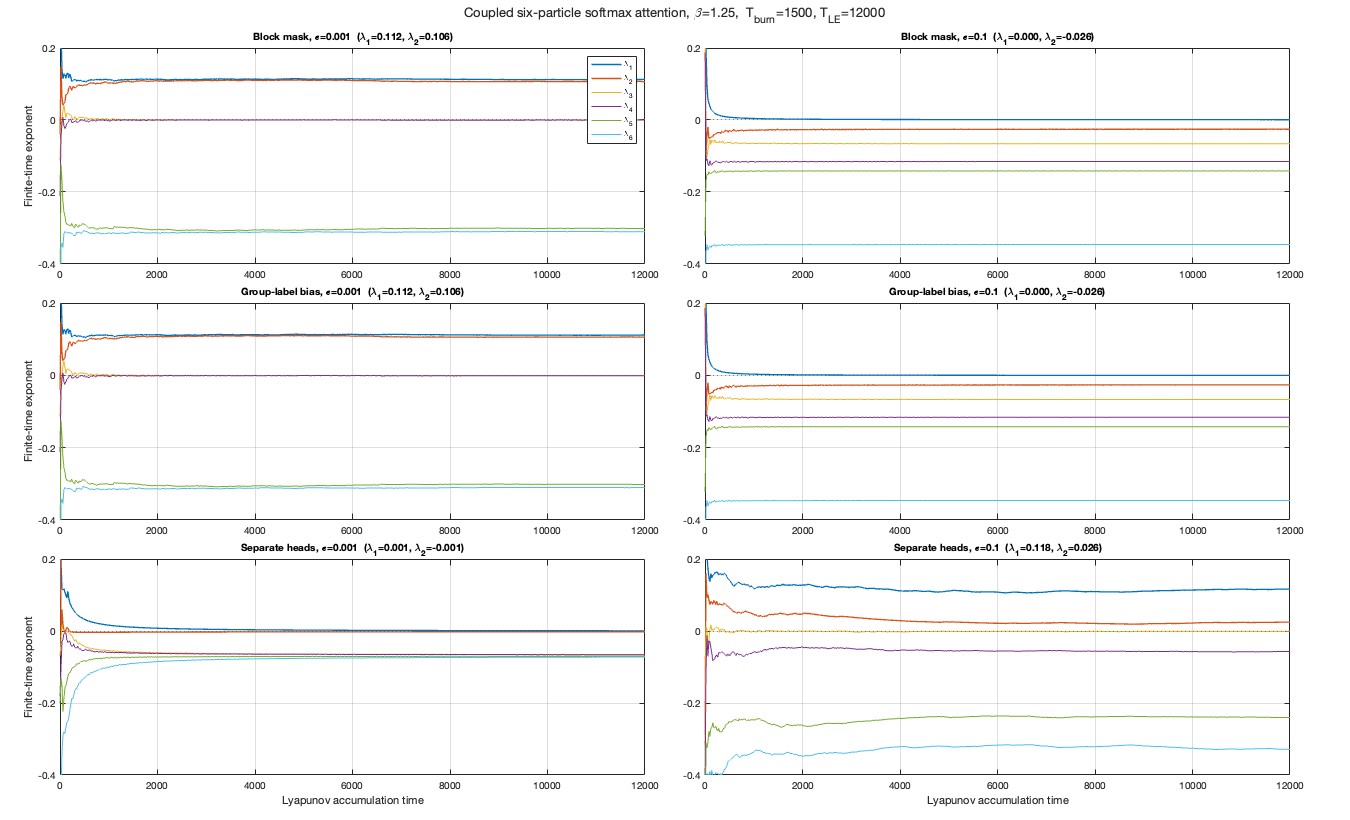}
    \caption{Finite-time Lyapunov exponents for two coupled
    three-particle systems on $\mathbb T^6$ at $\beta=1.25$.
    Columns correspond to $\varepsilon=0.001$ (left) and
    $\varepsilon=0.1$ (right); rows correspond to block-masked
    attention (top), fixed group labels (middle), and separately
    normalized internal and cross-group attention heads (bottom).
    The top two rows coincide exactly because the group-label bias
    $\gamma=-\log\varepsilon$ produces the same vector field as
    block-masked attention. For $\varepsilon=0.001$, the block-masked
    system retains two positive exponents, while the separate-head run
    displayed here has leading exponents near zero. At $\varepsilon=0.1$,
    the block-masked system has a spectrum consistent with an
    attracting periodic orbit, whereas the separate-head system has
    two positive exponents. Exponents are computed after discarding a
    transient of length $1500$ and are accumulated over $12000$ time
    units.}
    \label{fig:sixParticles}
\end{figure}

Next, we ask whether the chaotic dynamics found in the three-particle
model can survive when two such systems are coupled.  This is the simplest
nontrivial test of the constructions introduced in
Section~\ref{sec:bigTransformers}.  Consider six
particles, and partition them into two groups
 $G_1=\{1,2,3\}$ and $G_2=\{4,5,6\}$ with all particles
constrained to the ellipse determined by $W=\operatorname{diag}(2,1)$.
The phase space is $\mathcal E_W^6\cong\mathbb T^6$.

In each subgroup we choose effective matrices
\[
 A=\begin{pmatrix}-0.05&-2.30\\ 2.30&-0.05\end{pmatrix},
 \qquad
 B=\begin{pmatrix}-1.20&0.40\\ -0.80&-2.40\end{pmatrix},
 \qquad \beta=1.25.
\]
For these parameters, the uncoupled three-particle system exhibits
numerical evidence of a chaotic attractor, with Lyapunov spectrum
approximately $(0.11,0,-0.31)$.  At $\varepsilon=0$ the six-dimensional
system is the direct product of two copies of this flow.  In particular,
its Lyapunov spectrum contains two positive exponents and, up to
numerical error, two zero exponents.  Our initial configurations for
the coupled experiments are obtained by taking distinct points along
the previously computed chaotic three-particle trajectory in each group.
We compare the dynamics of the coupled six-dimensional system
with the three implementations (corresponding to two distinct
coupled vector fields) considered in Section~\ref{sec:bigTransformers}.  

For block-masked attention, the
cross-group entries of the softmax numerator are multiplied by
$\varepsilon$, while the within-group entries retain weight one; the
entire row is then normalized as in~\eqref{eq:big-softmask}.
In the fixed-group-label construction, the scores receive the bias
$\gamma\mathbf 1_{\{g(i)=g(j)\}}$.  The relation
$\gamma=-\log\varepsilon$ makes the resulting vector field
\emph{exactly identical} to the block-masked field for every
$\varepsilon>0$.  Hence these two constructions give the same
numerical experiment, although they represent different architectural
implementations.  

Finally, the two-head construction of Equation 
\eqref{eq:big-twoheads} combines a separately normalized internal head
with $\varepsilon$ times a separately normalized cross-group head.
For this first comparison, the two heads use the same matrices $A,B$
and inverse temperature $\beta$; thus the only new parameter is the
coupling strength $\varepsilon$.

We consider $\varepsilon=10^{-3}$ and $\varepsilon=10^{-1}$.
The angular equations and their analytic variational equations are
integrated using a fourth-order Runge--Kutta method with time step
$h=0.02$.  We discard an initial interval of length $1500$ and then
estimate all six Lyapunov exponents over a further interval of length
$12,000$, using QR orthogonalization at unit time intervals.
As a numerical consistency check, the sum of the estimated exponents
agrees with the time average of the divergence to approximately
$10^{-8}$ in the long integrations.  Figure~\ref{fig:sixParticles}
shows the convergence of the finite-time exponents for the six
combinations of implementation and coupling parameter.

For weak coupling, $\varepsilon=10^{-3}$, the block-masked system retains 
two positive Lyapunov exponents, with $(\lambda_1,\lambda_2)\approx(0.1125,0.1062)$. 
In contrast, the separate-head system has leading exponents approximately $(0.0015,-0.0013)$, 
providing no convincing numerical evidence of chaotic dynamics in this experiment. 
Thus, the two coupling architectures produce qualitatively different behavior even at this small parameter value.
At the larger coupling strength $\varepsilon=0.1$, the behavior of the
two models differs substantially, but in the other direction.  For the block-masked system
(equivalently, the group-label system), we find
$\lambda_1\approx 0.00031$ and $\lambda_2\approx-0.02614$,
with the remaining exponents also negative.  This spectrum is
consistent with convergence to an attracting periodic orbit rather
than a chaotic attractor.  Indeed, a subsequent numerical return
calculation gives a candidate period $T\approx64.27594$, with return
error of order $10^{-6}$. By contrast, the separate-head system has
a finite-time spectrum with two positive exponents, approximately
$(0.1176,0.0256)$ followed by an exponent near zero and three negative exponents.
This provides numerical evidence for hyperchaotic dynamics in the
six-dimensional flow. The contrast is visible in
Figure~\ref{fig:sixParticles}.

These experiments suggest that the observed chaotic dynamics depend
sensitively on the architecture used to couple the subsystems, and not
only on the presence of weak cross-group interactions.  Nevertheless,
$\varepsilon$ does not measure the same perturbation of the vector
field in the two models: the block-masked head renormalizes internal
and external weights together, whereas the cross-group head is
normalized separately and then added with coefficient $\varepsilon$.
Thus the comparison at equal $\varepsilon$ is a comparison of two
specific architectural parameterizations, not of perturbations with
equal norm.  Preliminary experiments at intermediate coupling values
also suggest that the observed dynamics can depend on the initial
condition, so possible coexistence of attracting sets deserves further
study. In what follows, the positive exponents observed for the
separate-head construction at $\varepsilon=0.1$ motivate experiments with larger
networks of coupled three-particle systems.

\subsection{Coupling in nine-particle systems}
\label{sec:higherDim}

We next consider a heterogeneous network in which a chaotic
three-particle system is coupled to two three-particle systems
possessing attracting periodic orbits. The resulting flow has nine
angular variables, with particle groups
$G_1=\{1,2,3\}$, $G_2=\{4,5,6\}$, and $G_3=\{7,8,9\}$.
All three subsystems use the same effective matrices
\[
A=\begin{pmatrix}-0.05&-2.30\\2.30&-0.05\end{pmatrix},
\qquad
B=\begin{pmatrix}-1.20&0.40\\-0.80&-2.40\end{pmatrix},
\]
but have distinct inverse-temperature parameters,
\[
(\beta_1,\beta_2,\beta_3)=(1.25,0.55,0.55).
\]
In isolation, the first subsystem exhibits a chaotic attracting set,
whereas each of the other two possesses the attracting nonsynchronized
periodic orbit discussed above. The uncoupled product therefore has
one positive Lyapunov exponent, associated with the chaotic component;
the remaining two subsystems are regular. This setup allows us to ask
whether coupling simply preserves the original chaotic direction or
creates additional expanding directions in the full network.

We compare block-masked attention with separately normalized internal
and cross-group attention heads, as defined in
Section~\ref{sec:bigTransformers}, fixing $\varepsilon=0.1$. The
parameters $A$, $B$, and $\beta_g$ used to evaluate each attention row
are those of its receiving group. For each architecture, we calculate
the Lyapunov spectrum from two distinct initial configurations.
The angular equations and their analytic variational equations are
integrated with fourth-order Runge--Kutta at step size $h=0.02$.
After discarding a transient of length $1500$, the exponents are
accumulated for $9000$ time units using periodic QR orthogonalization.
Finite-difference checks of the analytic Jacobian give relative errors
below $10^{-9}$. The results are summarized in
Table~\ref{tab:nineParticleLyapunov}.

\begin{table}[tbp]
\centering
\begin{tabular}{lcccc}
\hline
Coupling architecture & Initial condition & $\lambda_1$ & $\lambda_2$
    & $\#\{\lambda_j>0.01\}$\\
\hline
Block-masked & 1 & $\phantom{-}0.00004$ & $-0.09094$ & $0$\\
Block-masked & 2 & $-0.00010$ & $-0.09076$ & $0$\\
Separate heads & 1 & $\phantom{-}0.14384$ & $\phantom{-}0.09395$ & $3$\\
Separate heads & 2 & $\phantom{-}0.14544$ & $\phantom{-}0.09638$ & $3$\\
\hline
\end{tabular}
\caption{Nine-particle system consisting of one chaotic and two
attracting-periodic three-particle subsystems, coupled at
$\varepsilon=0.1$. The table reports the two largest finite-time
Lyapunov exponents and the number of exponents greater than $0.01$,
for two distinct initial configurations in each architecture.
The integration uses $h=0.02$, a discarded transient of length $1500$,
and $9000$ time units of Lyapunov accumulation.}
\label{tab:nineParticleLyapunov}
\end{table}

The contrast between the two architectures is pronounced. With
block-masked coupling, the largest exponent approaches zero, while
the second is approximately $-0.091$ for both initial conditions.
This is consistent with convergence to regular attracting dynamics,
and provides no numerical evidence of a surviving positive exponent.
By contrast, with separately normalized heads, the largest two
exponents stabilize near $0.14$ and $0.095$. Moreover, the computation
resolves \emph{three} exponents exceeding $0.01$ in both runs;
in particular, the third exponent exceeds this threshold. Thus
an uncoupled system with only one chaotic constituent produces, after
coupling, a nine-dimensional flow with three expanding directions.
The emergence of additional positive exponents is particularly
noteworthy because the two receiving subsystems have attracting
periodic dynamics when isolated.

We emphasize that positive Lyapunov exponents characterize expansion
in the coupled system; they do not, by themselves, identify which
constituent particles are responsible for the expanding directions.
Nevertheless, this experiment illustrates the central point of the
coupled-system study: knowledge of the dynamics of low-dimensional
building blocks provides a useful starting point, while the coupling
architecture can qualitatively reorganize their collective behavior.
The reported exponents are numerical estimates rather than validated
claims about an attractor.

\subsection{A $60$-particle example}
\label{sec:sixtyParticles}

We now investigate whether the chaotic dynamics observed in the
three-particle system persist in higher-dimensional systems, and
how interactions between chaotic subsystems influence the number
of positive Lyapunov exponents.

We consider a system of $N=60$ particles, organized into
$M=20$ groups of three particles each. Within each group,
the dynamics are determined by the same constant coefficient
matrices and inverse temperature $\beta=1.25$ used in the
three-particle chaotic example discussed above. The groups
are coupled through a multi-head attention mechanism consisting
of two separately normalized contributions: an internal
attention head, involving particles within the same group,
and an external attention head, involving particles in the
other groups. The external contribution is multiplied by
a coupling parameter $\varepsilon$.

When $\varepsilon=0$, the system reduces to a direct product
of twenty independent three-particle systems. Since each
particle is constrained to an ellipse, the resulting phase
space is $60$-dimensional. In the coupled experiments we
set $\varepsilon=0.1$.

The equations and their variational equations were integrated
numerically using a fourth-order Runge--Kutta method with
step size $h=0.0125$. An initial transient of length
$T_{\mathrm{tr}}=2000$ was discarded, after which the
largest $k=40$ Lyapunov exponents were approximated over
an additional integration interval of length $T=16000$.
We repeated the experiment for two independently generated
initial conditions. As a consistency check, the analytic
Jacobian was compared with finite-difference approximations,
giving relative errors between $6.3\times 10^{-10}$ and
$7.6\times 10^{-10}$ in the coupled experiments.

Table~\ref{tab:sixty_particle_lyapunov} summarizes the
results. In both experiments, the largest Lyapunov exponent
is approximately $0.13$, while the second largest is
approximately $0.12$. More significantly, both experiments
produce $26$ Lyapunov exponents exceeding $0.01$.
The estimates of the leading exponents remain relatively
stable during the final portion of the integration, and
the number of exponents above this threshold is unchanged
throughout most of the computation.

For comparison, we also compute the Lyapunov exponents of
an isolated three-particle subsystem. This calculation
gives a largest exponent of approximately $0.112$, a
second exponent close to zero, and one exponent exceeding
$0.01$. Consequently, the corresponding uncoupled product of twenty
such subsystems would have twenty exponents exceeding $0.01$,
counting multiplicity, at the level of these numerical estimates.

\begin{table}[ht]
\centering
\small
\begin{tabular}{lcccccc}
\hline
System & $N$ & $\varepsilon$ & $k$
& $\lambda_1$ & $\lambda_2$ & $n_{+}$ \\
\hline
Coupled, seed 0
& 60 & 0.1 & 40 & 0.1301503 & 0.1231389 & 26 \\
Coupled, seed 1
& 60 & 0.1 & 40 & 0.1297627 & 0.1221250 & 26 \\
Uncoupled, single group
& 3 & 0 & 3 & 0.1118411 & $-0.0000194$ & 1 \\
\hline
\end{tabular}
\caption{Numerically computed Lyapunov exponents for the
$60$-particle coupled system and the isolated three-particle
baseline. Here $N$ is the number of particles, $k$ is the
number of computed Lyapunov exponents, and $n_{+}$ denotes
the number of computed exponents greater than $0.01$.
All computations use a transient of length $2000$ and
a Lyapunov averaging interval of length $16000$.}
\label{tab:sixty_particle_lyapunov}
\end{table}

The comparison is particularly interesting because the
coupled system has six more numerically estimated Lyapunov
exponents exceeding $0.01$ than the corresponding uncoupled
product. Thus, the coupling not only preserves numerical evidence
of chaotic behavior but appears to introduce additional expanding
directions. In particular, these experiments provide numerical
evidence of high-dimensional chaotic dynamics arising through
interactions between low-dimensional chaotic subsystems. This
comparison does not establish a general scaling law as the number
of groups increases.

We emphasize that only the largest $40$ of the $60$
Lyapunov exponents were computed. Since these are the largest
exponents and only $26$ exceed $0.01$, the omitted exponents
cannot exceed this threshold. Thus truncation does not affect
the reported count of exponents above $0.01$; the number of
strictly positive exponents could be larger because some may
lie between zero and $0.01$. The estimates are finite-time
numerical quantities. Nevertheless, the agreement between the
two independent experiments provides strong numerical evidence
that the coupled system supports many expanding directions.
These Lyapunov exponents describe the full tangent dynamics and
do not, by themselves, determine whether the trajectory lies
near an invariant synchronization subspace arising from the
permutation symmetry among identical groups.

\section{Conclusions}
\label{sec:conclusions}

The interpretation of transformers as interacting particle
systems opens the door to studying their behavior using the
geometric and analytic tools of dynamical systems theory.
In particular, the gradient structures identified in certain
simplified attention models provide invaluable insight into
the emergence of clustering and metastability. When a
gradient structure is available, the associated energy
landscape constrains the possible dynamics and provides
a natural framework for understanding long transients near
saddle equilibria and their stable manifolds. This point
of view has already led to important advances in the
mathematical theory of transformers.

However, gradient structure is not a prerequisite for
clustering or metastability. More generally, trajectories
may approach the stable manifolds of saddle-type invariant
sets, remain in their neighborhoods for extended periods,
and eventually depart along weakly unstable directions.
Such behavior does not require an underlying energy
function or a gradient-like flow. Moreover, the invariant
sets organizing these dynamics need not be equilibria.
Periodic orbits, invariant tori, and more complicated
saddle-type invariant sets may also play a role.

The analytic results in this paper show that, when the 
assumptions of Theorem \ref{thm:intro-gradient} are violated,
the self-attention system can exhibit non-trivial recurrent 
dynamics.  The numerical experiments show that even relatively
simple autonomous transformer models, with particles
constrained to ellipses and constant coefficient matrices,
can already exhibit rich dynamics. In particular, we
observe numerical evidence of chaos in low-dimensional
examples, and the persistence and proliferation of
expanding directions in larger coupled systems.
Our $60$-particle experiments provide
evidence of high-dimensional chaos, with $26$
numerically computed Lyapunov exponents exceeding
$0.01$. These observations suggest that complicated
recurrent behavior is not confined to low-dimensional
examples, and it may even become increasingly prevalent as
the number of interacting particles grows.

We speculate that complicated dynamics may
be a feature rather than a bug. After all,
a central purpose of the attention mechanism is to
produce context-dependent representations of information.
It is conceivable that the stretching and folding
associated with chaotic dynamics
aid in this task. Whether such dynamics
actually play a useful role in trained transformer
architectures is an intriguing question which the
present work does not address.

Many fundamental mathematical questions remain open.
We conjecture that metastable clustering in more general
attention systems may be organized by the stable manifolds
of saddle-type invariant sets, including periodic orbits,
invariant tori, and possibly chaotic invariant sets.
Identifying such structures and developing methods for
computing their stable and unstable manifolds, slow
invariant submanifolds, and associated invariant bundles
for self-attention systems 
are natural directions for future research. Understanding
how these structures depend on the attention parameters,
and how they interact through bifurcations, may provide
a more general geometric theory of clustering and
metastability beyond the gradient setting.

Our higher-dimensional experiments also suggest the
possibility of a systematic relationship between the
number of particles and the complexity of the dynamics.
In particular, we conjecture that suitable coupling
constructions can produce families of systems for which
the number of positive Lyapunov exponents grows with
the number of coupled subsystems. Establishing such
results, and understanding the persistence of chaotic
invariant sets under coupling, would provide a natural
extension of the present numerical investigations.

Finally, it is important to ask how much of this
dynamical structure survives in more realistic
transformer architectures. 
How does the inclusion of time 
varying effects in the query, key, and value matrices
change the geometric picture?
Can analogous invariant
sets and their stable and unstable manifolds be
identified in trained models? Can their geometry
help explain the organization and evolution of
token representations, and ultimately contribute
to the interpretability of these systems?

We regard the present work as an initial exploration of these questions. 
The mathematical theory of transformers is still in its early stages, 
but perhaps the relationship between attention and invariant dynamics 
itself deserves more attention.

\section{Acknowledgments} \label{sec:ack}

The author is deeply indebted to Yu Xiang
of Florida Atlantic University for introducing him to the 
particle system interpretation of LLM transformers studied in the present work, and 
for hours of indispensable conversation.  This paper would not exist without these
discussions.  Conversations with Juan Miranda and Kyle Loh are also 
gratefully acknowledged.  The author was partially supported by the National 
Science Foundation grant DMS - 2307987.

\section{Declaration of generative AI and AI-assisted technologies} \label{sec:AI}

Generative artificial intelligence tools, primarily ChatGPT (OpenAI), were used during every 
phase of the development of this work, including the generation, revision, and debugging of 
numerical code, assistance with mathematical calculations and technical discussions, and the 
drafting, editing, and organization of portions of the final manuscript.
Claude was also used for independent copy-edit review.

The use of AI was interactive and author-directed, involving many short chat mode 
exchanges in which specific mathematical, computational, or expository tasks were posed, 
the resulting outputs were examined, and subsequent redirections and revisions were requested. 
The research was not delegated to autonomous AI agents.

In particular, substantial pieces of numerical software were developed with AI assistance. 
The author tested and executed the resulting programs, performed the numerical experiments, 
examined the computed output, produced and interpreted the figures and other reported results. 
AI-generated text was similarly reviewed, edited, and incorporated by the author.

The author was responsible for the scientific questions, the development and selection of the 
mathematical models, the design and evaluation of the numerical experiments, and the interpretation 
of the results. All mathematical claims, numerical findings, and conclusions presented in the final 
manuscript are the responsibility of the human author, who takes full responsibility for the accuracy, 
originality, and integrity of the published work.


\end{document}